\documentclass[a4paper,11pt]{article}
\usepackage{graphicx}
\usepackage{a4wide}
\usepackage[english]{babel}
\usepackage{amsfonts}
\usepackage{amsmath}
\usepackage{amssymb}
\usepackage{dsfont}
\newtheorem{lemma}{Lemma}
\newtheorem{theorem}{Theorem}
\newtheorem{corollary}{Corollary}

\usepackage{float}

\newcommand {\E} {\mathbb{E}}

\DeclareMathOperator {\var}{Var_{[0,1]}}

\newcommand {\p} {\mathbb{P}}
\newcommand {\Z} {\mathbb{Z}}
\newcommand {\N} {\mathbb{N}}

\newcommand {\F} {\mathcal{F}}
\newcommand {\R} {\mathbb{R}}
\newcommand {\Q} {\mathbb{Q}}

\newcommand {\ve} {\varepsilon}

\makeatletter\def\blfootnote{\xdef\@thefnmark{}\@footnotetext}\makeatother
\allowdisplaybreaks
\title{\bf Quantitative uniform distribution results for geometric progressions}
\author{Christoph Aistleitner\footnote{Graz University of Technology,
Institute of Mathematics A, Steyrergasse 30, 8010 Graz, Austria. \mbox{e-mail}:
\texttt{aistleitner@math.tugraz.at}. The author is supported by a Schr\"odinger scholarship of the Austrian Research
Foundation FWF, as well as by the FWF-Project P24302.}}
\begin{document}

\date{}
\maketitle

\blfootnote{{\bf keywords:} geometric progression, discrepancy, law of the iterated logarithm,
Koksma's theorem, central limit theorem}

\begin{abstract}
By a classical theorem of Koksma the sequence of fractional parts $(\{x^n\})_{n
\geq 1}$ is uniformly distributed for almost all values of $x>1$. In the present
paper we obtain an exact quantitative version of Koksma's theorem, by
calculating the precise asymptotic order of the discrepancy of $(\{\xi
x^{s_n}\})_{n \geq 1}$ for typical values of $x$ (in the sense of Lebesgue
measure). Here $\xi>0$ is an arbitrary constant, and $(s_n)_{n \geq 1}$ can be
any increasing sequence of positive integers.
\end{abstract}

\section{Introduction and statement of results}

A sequence $(x_n)_{n \geq 1}$ of real numbers from the unit interval is called
\emph{uniformly distributed modulo 1} (u.d. mod 1) if for any $0 \leq a < b \leq
1$
\begin{equation} \label{konv}
\frac{1}{N} \sum_{n=1}^N \mathds{1}_{[a,b)} (x_n) \to b-a \qquad \textrm{as $N
\to \infty$}.
\end{equation}
In other words, a sequence is u.d. mod 1 if the relative number of elements of
the sequence contained in an interval $[a,b) \subset [0,1)$ always converges
to the length (or Lebesgue measure) of this interval. Here the length of such an
interval can be interpreted as the expected value for the relative number of
elements of a random sequence contained in it, and with regard to the
Glivenko--Cantelli theorem a uniformly distributed sequences can be considered as
a sequence showing random behavior. There exist many sequences which are u.d.
mod 1, for example the sequence $(\{n x\})_{n \geq 1}$ whenever $x \not\in \Q$
(here, and in the sequel, $\{ \cdot \}$ denotes the fractional part).\\

The speed of convergence in~\eqref{konv} is measured by the \emph{discrepancy}
and the \emph{star-discrepancy} of the sequence $(x_n)_{n \geq 1}$. For a finite
sequence $(x_1, \dots, x_N)$ the (extremal) discrepancy $D_N$ and the
star-discrepancy $D_N^*$ are defined as
$$
D_N(x_1, \dots, x_N) = \sup_{0 \leq a < b \leq 1} \left| \frac{1}{N}
\sum_{n=1}^N \mathds{1}_{[a,b)} (x_n) - (b-a) \right|
$$
and
$$
D_N^* (x_1, \dots, x_N) \sup_{0 < a \leq 1} \left| \frac{1}{N} \sum_{n=1}^N
\mathds{1}_{[0,a)} (x_n) - a \right|.
$$
For simplicity, we will write $D_N(x_n)$ and $D_N^*(x_n)$ for the discrepancy
resp. star-discrepancy of the first $N$ elements of a (finite or infinite)
sequence. For an introduction to the theory of uniform distribution modulo 1 and
discrepancy theory the reader is referred to the monographs~\cite{drti,kn}.\\

By a remarkable result of Weyl~\cite{weyl} for any sequence of distinct integers
$(s_n)_{n \geq 1}$ the sequence $(\{s_n x\})_{n \geq 1}$ is u.d. mod 1 for
almost all $x$ (in the sense of Lebesgue measure). This is equivalent to the
fact that
$$
D_N (\{s_n x\}) \to 0 \qquad \textrm{as $N \to \infty$} \qquad \textrm{for
almost all $x$}.
$$
Precise results are only known in a few special cases. For example, when
$s_n=n,~n \geq 1$, we have
\begin{equation} \label{kesten}
\frac{N D_N(\{n x\})}{\log N \log \log N} \to \frac{2}{\pi^2} \qquad \textrm{in
measure}
\end{equation}
due to Kesten~\cite{kesten} (see also~\cite{schoi}). Exact results of this type
are possible since there is an intimate connection between the discrepancy of
$(\{n x\})_{n \geq 1}$ and the continued fraction expansion of $x$. The second
class of sequences for which precise metric results are known are sequences
satisfying the Hadamard gap condition
$$
\frac{s_{n+1}}{s_n} \geq q > 1, \qquad n \geq 1.
$$
In this case Philipp~\cite{plt} proved the bounded law of the iterated logarithm
(LIL)
\begin{equation} \label{phillil}
\frac{1}{4} \leq \frac{\sqrt{N} D_N(\{s_n x\})}{\sqrt{\log \log N}} \leq C_q
\qquad \textup{a.e.},
\end{equation}
where $C_q$ depends only on the growth factor $q$ (the lower bound follows from
an older result of Erd\H os and G{\'a}l~\cite{egal} and Koksma's inequality). For sub-exponentially growing $(s_n)_{n \geq 1}$ the LIL~\eqref{phillil} generally fails, unless $(s_n)_{n \geq 1}$ satisfies some strong number-theoretic conditions (see for example~\cite{aiill,bpte}). 
For sequences of the special form $s_n = \beta^n,~n \geq 1$, for some $\beta>1$,
Fukuyama~\cite{ft} recently proved the precise LIL
$$
\frac{\sqrt{N} D_N(\{\beta^n x\})}{\sqrt{\log \log N}} = \sigma_\beta \qquad
\textup{a.e.},
$$
where $\sigma_\beta$ is a constant depending on the number-theoretic properties
of $\beta$ in a very complicated and interesting way. In particular
$$
\frac{\sqrt{N} D_N(\{2^n x\})}{\sqrt{\log \log N}} = \frac{2\sqrt{21}}{9} \qquad
\textup{a.e.},
$$
and
\begin{equation} \label{fukuex}
\frac{\sqrt{N} D_N(\{\beta^n x\})}{\sqrt{\log \log N}} = \frac{1}{\sqrt{2}}
\qquad \textup{a.e.}
\end{equation}
if $\beta$ is a number for which $\beta^r \not\in \Q$ for all $r \geq 1$. These results should be
compared to the Chung--Smirnov law of the iterated logarithm for independent,
identically $[0,1]$-uniformly distributed random variables $(X_n)_{n \geq 1}$,
which states that
\begin{equation*} \label{chungs}
\frac{\sqrt{N} D_N(X_n)}{\sqrt{\log \log N}} = \frac{1}{\sqrt{2}} \qquad
\textup{a.s.}
\end{equation*}
In this specific form the law of the iterated logarithm for the discrepancy (in the language of probability theory: for the Kolmogorov--Smirnov statistic) of $(X_n)_{n \geq 1}$
is due to Chung~\cite{chung} and Cassels~\cite{casselsLIL}; for a general
formulation, see e.g.~\cite[p.~504]{show}. Recall that a number $x$ is a normal
number in base $\beta$ if and only if $D_N(\{\beta^n x\}) \to 0$ as $N \to
\infty$. Consequently Fukuyama's result is a precise quantitative version of
Borel's well-known theorem that almost all numbers are normal~\cite{borel}.\\

To the best of my knowledge the two mentioned classes of sequences (arithmetic
progressions and lacunary sequences) are essentially the only two classes of
parametric sequences for which the typical (in the sense of Lebesgue measure)
asymptotic order of the discrepancy  is precisely known. For general sequences
$(s_n)_{n \geq 1}$ of distinct integers we only have the upper bounds
$$
D_N(\{s_n x\}) = \mathcal{O} \left( \frac{(\log N)^{5/2+\ve}}{\sqrt{N}} \right) \qquad
\textup{a.e.}
$$
(Erd\H os and Koksma~\cite{erko2}) and
\begin{equation} \label{bak}
D_N(\{s_n x\}) = \mathcal{O} \left( \frac{(\log N)^{3/2+\ve}}{\sqrt{N}} \right) \qquad
\textup{a.e.} \qquad \textrm{if $(s_n)_{n \geq 1}$ is increasing}
\end{equation}
(Baker~\cite{baker2}). It is known that the exponent of the logarithmic term in
\eqref{bak} can in general not be reduced below 1/2 (Berkes and
Philipp~\cite{bp}), but as~\eqref{kesten} shows for a specific sequence
$(s_n)_{n \geq 1}$ the typical speed of convergence of $D_N(\{s_n x\})$ can
differ from~\eqref{bak} significantly. More details on metric discrepancy theory
can be found in the book of Harman~\cite{harman} and in the survey paper~\cite{asus}.\\

In 1935, Koksma~\cite{kok} proved a very general result in uniform distribution theory, which as a special case contains the fact that for any $\xi>0$ and any sequence $(s_n)_{n \geq 1}$ of distinct positive integers the sequence $(\{\xi x^{s_n}\})_{n \geq 1}$ is u.d. mod 1 for almost all $x>1$. In particular, geometric progressions $(\{x^n\})_{n \geq 1}$ are u.d. mod 1 for almost all $x>1$. Erd\H os and Koksma~\cite{erko1} proved that the asymptotic order of the discrepancy of $(\{\xi x^{s_n}\})_{n \geq 1}$, in the case of increasing $(s_n)_{n \geq 1}$, satisfies
\begin{equation} \label{erko}
D_N(\{\xi x^{s_n}\}) = \mathcal{O} \left( \frac{(\log N)^{3/2} (\log \log N)^{1/2+\ve}}{\sqrt{N}} \right) \quad \textrm{as $N \to \infty$} \quad \textrm{for almost all $x>1$}.
\end{equation}
In 1950 this was improved by Cassels~\cite{cassels}, who obtained
\begin{equation} \label{cassels}
D_N(\{\xi x^{s_n}\}) = \mathcal{O} \left( \frac{\log N (\log \log N)^{3/2+\ve}}{\sqrt{N}}\right) \quad \textrm{as $N \to \infty$} \quad \textrm{for almost all $x>1$}.
\end{equation}
Since then, no further improvements of~\eqref{cassels} have been made. On the other hand, as far as I know, no asymptotic \emph{lower} bounds for $D_N(\{\xi x^{s_n}\})$ or $D_N(\{x^{n}\})$ for typical values of $x$ (in the sense of Lebesgue measure) have ever been proved.\\

Concerning the asymptotic distribution of $(\{x^{s_n}\})_{n \geq 1}$, it should be mentioned that Niederreiter and Tichy~\cite{nt2} proved that this sequence is completely uniformly distributed\footnote{A sequence $(x_n)_{n \geq 1}$ is called completely uniformly distributed modulo 1, if for any $s \geq 1$ the $s$-dimensional sequence $((x_n,\dots,x_{n+s-1}))_{n \geq 1}$ is uniformly distributed mod 1 in $[0,1]^s$. See~\cite{drti,kn} for details.} modulo~1 for almost all $x>1$, by this means solving a problem posed by Knuth~\cite{knuth}, who suggested complete uniform distribution as a criterion for pseudorandomness of deterministic sequences. For the case $s_n=n$ this had already been shown by Franklin~\cite{franklin}. For more details see~\cite{nt2}, as well as~\cite{nt1,nt3} for quantitative results containing discrepancy estimates. Koksma's theorem has been generalized to many other cases, including complex numbers (LeVeque~\cite{lev}), quaternions (Tichy~\cite{tiq} and
Nowak~\cite{noq}) and matrices (Nowak and Tichy~\cite{noti1,noti2}). In many cases, quantitative results similar to~\eqref{erko} and~\eqref{cassels} exist.\\

Deciding whether $(\{ x^n\})_{n \geq 1}$ or $(\{\xi x^{s_n}\})_{n \geq 1}$ is u.d. mod 1 for a specific value of $x$ is a notoriously difficult problem, and only few partial results are known. For a comprehensive survey, see~\cite{bug}. A famous open problem is whether the sequence $(\{(3/2)^n\})_{n \geq 1}$ is u.d. mod 1, but in fact we do not even know if $\limsup_{n \to \infty} \{(3/2)^k\} - \liminf_{n \to \infty} \{(3/2)^k\} \geq 1/2$. Mahler~\cite{mahler} asked whether there exists a number $\xi>0$ for which $\{\xi (3/2)^n \} \in [0,1/2]$ for all $n \geq 1$. These problems are connected with other difficult mathematical problems, such as Waring's problem (see~\cite{ben}) and the $3x+1$ problem (see~\cite{fla}). For an overview see~\cite{strauch}; recent contributions are e.g. due to Akiyama, Frougny and Sakarovitch~\cite{akie,aki}, Dubickas~\cite{dub1,dub2} and Kaneko~\cite{kaneko1,kaneko2}.\\

It is known that in some cases $(\{\xi x^n\})_{n \geq 1}$ fails to be uniformly distributed. This is for example the case when $x$ is a Pisot--Vijayaraghavan number and $\xi$ is an algebraic integer in the field of $x$; in this case for the distance to the nearest integer $[\cdot ]$ we have $[\xi x^n] \to 0$ at an exponential rate. A longstanding problem of Hardy asks whether there exists a transcendental $x>1$ for which there is some $\xi>0$ such that $[\xi x^n] \to 0$. A recent result of Bugeaud and Moshchevitin~\cite{bugm} in this context is the following: Modifying a probabilistic method of Peres and Schlag~\cite{pst}, they proved that there exist arbitrarily small numbers $\ve>0$ such that $[(1+\ve)^n] > 2^{-17} \ve |\log \ve|^{-1}$ for all $n \geq 1$. For more details on these problems the reader is again referred to~\cite{bug}.\\

The purpose of the present paper is to prove precise metric results for the
discrepancy of sequences of the form $(\xi x^{s_n})_{n \geq 1}$, where $\xi>0$
is a arbitrary (fixed) number and $(s_n)_{n \geq 1}$ is an arbitrary increasing sequence of
positive integers. As a special case we obtain precise metric results
for geometric progressions $(\xi x^{n})_{n \geq 1}$, and as a byproduct of the
proof of our main theorem we also obtain a central limit theorem for $(f(\xi
x^{s_n}))_{n \geq 1}$.\\

As mentioned before, Theorem~\ref{th1} does not only improve the known metric upper bounds for the discrepancy of sequences of the form $(\{ \xi x^{s_n} \})_{n \geq 1}$, but also provides the very first metric lower bounds for the discrepancy of such sequences. In particular Theorem~\ref{th1} solves a problem of V.I.~Arnold, who in one of his final papers formulated the conjecture that the discrepancy $D_N$ of $(\{ \xi x^{n} \})_{n \geq 1}$ is \emph{not} of order $o(N^{-1/2})$ for almost all $x>1$ (see~\cite[p.~36]{arnold}). 

\begin{theorem} \label{th1}
For any strictly increasing sequence $(s_n)_{n \geq 1}$ of positive integers and any number
$\xi>0$ we have for almost all $x>1$
$$
\limsup_{N \to \infty} \frac{\sqrt{N} D_N(\{ \xi x^{s_n} \})}{\sqrt{\log \log
N}} = \limsup_{N \to \infty} \frac{\sqrt{N} D_N^*(\{ \xi x^{s_n} \})}{\sqrt{\log
\log N}} = \frac{1}{\sqrt{2}}.
$$
\end{theorem}

As an immediate consequence of Theorem~\ref{th1} we obtain the following
corollary for geometric progressions of the form $(\xi x^n)_{n \geq 1}$. 
\begin{corollary} \label{co1}
For any $\xi>0$ we have for almost all $x>1$ 
$$
\limsup_{N \to \infty} \frac{\sqrt{N} D_N(\{ \xi x^n \})}{\sqrt{\log \log N}} =
\limsup_{N \to \infty} \frac{\sqrt{N} D_N^*(\{ \xi x^n \})}{\sqrt{\log \log N}}
= \frac{1}{\sqrt{2}}.
$$
\end{corollary}

As a byproduct of our proof of Theorem~\ref{th1} we also get the
following central limit theorem.
\begin{theorem} \label{th2}
Let $f$ be a function satisfying 
\begin{equation} \label{f}
f(x) = f(x+1), \qquad \int_0^1 f(x)~dx = 0, \qquad \var f \leq 2.
\end{equation} 
Then for any sequence $(s_n)_{n \geq 1}$ of distinct positive integers, any number
$\xi>0$ and any nonempty interval $[A,B] \subset (1,\infty)$ we have
$$
\p \left( x \in [A,B]:~\frac{1}{\sqrt{N}} \sum_{n=1}^N f(\xi x^{s_n}) \leq t\right) \to \Phi(t) \qquad \textrm{as $N \to \infty$}.
$$
Here $\p$ denotes the normalized Lebesgue measure on $[A,B]$, and $\Phi$ denotes the standard normal distribution function. The convergence is uniform in $t \in \R$.
\end{theorem}

I am not certain if these results are expected or surprising. Of course, it is
reasonable to imagine that the functions $f(\xi x^m\}$ and $f(\xi x^n)$ are ``almost
independent'' if the difference between $m$ and $n$ is large, and that therefore
the system $(f(\xi x^{s_n}))_{n \geq 1}$ and the discrepancy $D_N(\{\xi
x^{s_n}\})$ should show ``almost'' the same behavior as in the case of an i.i.d.
random sequence. For example, Beck~\cite[p.~55]{beck} writes that Koksma's
theorem on the uniform distribution of $(\{x^n\})_{n \geq 1}$ for a.e. $x$ ``was extended later to
more delicate results such as the law of the iterated logarithm and the central
limit theorem'', although such results have not been proved so far; apparently
Beck was convinced that they must be true.\\

However, comparing the case of
sequences of the form $(\xi x^{s_n})_{n \geq 1}$ to the somewhat similar case of lacunary
sequences, one sees that it is by no means clear that the (precise) LIL and CLT
have to hold for geometric progressions. In the case of lacunary sequences $(s_n
x)_{n \geq 1}$, the value of the limsup in the LIL for the discrepancy depends
on the precise number-theoretic properties of $(s_n)_{n \geq 1}$ in a very
complicated way, and can even be non-constant (see~\cite{airr1,airr2,airr3}). Furthermore, the asymptotic behavior of lacunary sequences can change significantly after a permutation of its terms, see~\cite{fuku1,fuku2}. 
Similarly, the CLT for lacunary sequences $(s_n)_{n \geq 1}$ is only true if the
sequence satisfies certain number-theoretic conditions, and the limit
distribution of $N^{-1/2} \sum_{n=1}^N f(s_n x)$ can fail to be Gaussian
(see~\cite{a1}).\\

Lacunary sequences and geometric progressions are essentially
of the same order of growth, so it could also be imagined that additional
number-theoretic conditions (like the Diophantine conditions in the case of
lacunary sequences) would be necessary to obtain the precise LIL and CLT for geometric progressions.
However, no such additional conditions are necessary, and apparently this is due to the
fact that sequences of the form $(\{\xi x^{s_n}\})$  by construction necessarily
have a more inhomogeneous structure than lacunary sequences, which for example
in the case $(\{2^n x\})_{n \geq 1}$ can have a very strong periodic and
homogeneous structure with respect to both $x$ and $n$. For the relation between the metric discrepancy results for geometric progressions in this paper and similar metric discrepancy results for lacunary series see also the addendum at the end of this paper, which I owe Katusi Fukuyama.\\

The proof of the main theorem of this paper is based on methods which were
developed for lacunary function systems. However, there are several
major differences to the case of geometric progressions, which made it necessary
to develop a new machinery. The two most significant differences are:
\begin{itemize}
 \item Lacunary systems $(f(s_n x))_{n \geq 1}$ have a direct connection with
Fourier analysis, and can be expanded into a Fourier series in a very simple and
natural way. This makes it possible to reduce the calculation of $L^p$-norms or
exponential norms to counting the number of solutions of Diophantine equations,
by utilizing the orthogonality of the trigonometric system. In the case of geometric progressions this is not possible, and instead of orthogonality properties we have
to use the fact that a function $f(\xi x^n)$ is highly oscillatory in comparison
with $f(\xi x^m)$ if $n \gg m$. While in the lacunary case the orthogonality of
the trigonometric system guarantees that in calculating integrals most of the
mixed factors vanish, we have to use the van der Corput inequality (see below)
instead and take care of a huge number of small quantities.
 \item For any $f$ satisfying~\eqref{f} and any integer $n$ the function $f(n
x)$ is periodic with period $n$, which means that
the global problem of considering all possible values of $x$ can often be
reduced to considering $x$ only ``locally'', and all values of $x$ can be treated in the same way. In the present case we have
functions of the form $f(\xi x^n)$, which do not posses this homogeneous
structure. On the contrary, the speed of oscillation of $f(\xi x^n)$ increases
as $x$ increases, which for example makes the martingale approximation in
Section~\ref{sectmar} much more complicated than in the lacunary case.
\end{itemize}

The rest of this paper is organized as follows:\\

\begin{itemize}
\item In Section~\ref{pre} we formulate several auxiliary results which will be
necessary for the proofs. In particular this includes the statement of the van
der Corput lemma, which is a crucial ingredient in our proof.
\item In Section~\ref{expi} we prove a large deviations bound for
$\sum_{n=M+1}^{M+N} f(\xi x^{s_n})$, which is a consequence of an exponential
inequality in the spirit of Takahashi~\cite{taka} and Philipp~\cite{plt}.
\item In Section~\ref{maxi} we prove a maximal version of the large deviations
inequality from Section~\ref{expi}. For the proof of this maximal inequality we
use a dyadic decomposition of the index set.
\item In Section~\ref{smallr} we prove a bounded law of the iterated logarithm
for functions which are the remainder of a Fourier series of a function
satisfying~\eqref{f}. Since the contribution of the remainder function of the
$d$-th partial sum of the Fourier series is small, it is sufficient to prove the
exact LIL for trigonometric polynomials instead of general functions $f$.
\item In Section~\ref{smalln} we prove a bounded law of the iterated logarithm
for a modified discrepancy, which considers only ``small'' subintervals of
$[0,1]$. We show that the contribution of these small intervals is small,
and that the proof of Theorem~\ref{th1} can be reduced to proving the exact LIL
for a single function $f$ instead of a supremum over uncountable many indicator
functions.
\item In Section~\ref{sectmar} and Section~\ref{sectp} we prove the exact LIL
for trigonometric polynomials. The proof uses an approximation by martingale
differences, which has been developed by Berkes and, independently, Philipp and Stout. The main
ingredient in the proof is a martingale version of the Skorokhod representation theorem due to Strassen.
\item In Section~\ref{sectf} the precise LIL for functions satisfying~\eqref{f}
is obtained as a consequence of the results from Section~\ref{smallr} and
Section~\ref{sectp}.
\item Finally, in Section~\ref{sectth1} we give the proof of Theorem~\ref{th1}.
In Section~\ref{sectth2} we show how the proof of Theorem~\ref{th2} can be obtained without
much additional effort as a byproduct of the proof of Theorem~\ref{th1}.
\end{itemize}

\section{Preliminaries} \label{pre}

We will assume throughout the rest of the paper that the number $\xi>0$ is
fixed. Furthermore, it is sufficient to prove that Theorem~\ref{th1} holds for
almost all $x \in [A,B]$, where $[A,B] \subset (1,\infty)$ is an arbitrary
interval. Throughout the rest of the paper, the numbers $A,B$ satisfying
$1 < A < B$ will be fixed. We will write $c$ for positive numbers, not always
the same, which may \emph{only} depend on $\xi$ and $A,B$, but not on $N,n,f,d$
or anything else (unless stated otherwise at the beginning of the respective
section). In the same sense we will use the symbols ``$\ll$'' and ``$\gg$''. For simplicity of
writing we will assume that $B-A=1$, which means that the interval $[A,B]$, equipped with Borel sets and Lebesgue measure, is a probability space. We will write $\p$ for the Lebesgue measure on $[A,B]$, and $\E$ for the expected value with respect to this measure.\\

Throughout the rest of this paper, we will write $\exp(x)$ for $e^x$.
Furthermore, $\log x$ denotes the natural logarithm, and should be interpreted as $\max\{1,\log x\}$. We set
$$
\|f\| = \left(\int_A^B (f(x))^2~dx \right)^{1/2}
$$
and 
\begin{equation} \label{I}
\mathbf{I}_{[a,b)} (x) = \mathds{1}_{[a,b)} (x) - (b-a).
\end{equation} 
Then for any $0 \leq a < b \leq 1$ the function $\mathbf{I}_{[a,b)}$ satisfies
\eqref{f}.

\begin{lemma}[{\cite[p.~48]{zt}}] \label{ztlemma}
Let $f$ be a function satisfying~\eqref{f}, and write 
$$
f(x) \sim \sum_{j=1}^\infty \left(a_j \cos 2 \pi x + b_j \sin 2 \pi x\right)
$$ 
for its Fourier series. Then
$$
|a_j| \leq \frac{1}{j}, \qquad |b_j|\leq \frac{1}{j}, \qquad \textrm{for}\quad j
\geq 1.
$$
\end{lemma}

The following Lemma~\ref{vclemma} is a special case of the van der Corput lemma.
It can be found e.g. in~\cite[Chapter 1, Section 1, Lemma 2.1]{kn}
or~\cite[Chapter VIII, Proposition 2]{stein}.
\begin{lemma} \label{vclemma}
Suppose that $\phi(x)$ is real-valued, that $|\phi'(x)| \geq \gamma$ for some
positive $\gamma$, and that $\phi'$ is monotonic for all $x \in (\alpha,
\beta)$. Then
$$
\left|\int_\alpha^\beta e^{2 \pi i \phi(x)}~dx\right| \leq \gamma^{-1}.
$$
\end{lemma}

Lemma~\ref{lemma3} and Lemma~\ref{lemma4} below follow directly from
Lemma~\ref{vclemma}.
\begin{lemma} \label{lemma3}
Let $n$ be a positive integer. Then for any subinterval $[\alpha,\beta]$ of
$[A,B]$ and any integer $j \geq 1$ we have
$$
\left| \int_\alpha^\beta  \cos (2 \pi j \xi x^n)~dx \right| \leq \frac{1}{j \xi n
\alpha^{n-1}}.
$$
\end{lemma}

\begin{lemma} \label{lemma4}
Let $m \neq n$ be positive integers. Then for any positive integers $j, k$ and any subinterval $[\alpha,\beta]$ of $[A,B]$ we have
$$
\left| \int_\alpha^\beta \cos (2 \pi \xi (j x^n + k x^m)) ~dx \right| \leq
\frac{1}{\xi \max\{m,n\} \alpha^{\max\{m,n\}-1}}.
$$
\end{lemma}

\begin{lemma} \label{lemma5}
Let $m < n$ be positive integers. Then for any positive integers $j,k$ and for
any $\eta>0$ there exist three disjoint intervals $I_1, I_2, I_3$ (depending on $j,k,m,n$) such that 
$$
\p \left([A,B] \backslash (I_1 \cup I_2 \cup I_3)\right) \leq 2 B \eta
$$
and such that for any interval $[\alpha,\beta]$ which is completely contained in one of the
intervals $I_1,~I_2$ or $I_3$ we have
$$
\left| \int_\alpha^\beta \cos (2 \pi \xi (j x^n - k x^m)) ~dx \right| \leq
\frac{1}{\xi \eta m \alpha^{m-1}}.
$$
\end{lemma}

\emph{Proof of Lemma~\ref{lemma5}:}~We want to use Lemma~\ref{vclemma} for
$\phi(x) = 2 \pi \xi(j x^n - k x^m)$. Obviously this is not directly possible, since it
might happen that $\phi'(x) = 2\pi\xi(jnx^{n-1} - kmx^{m-1}) =0$ for some $x$. We have
$$
\phi'(x)= \xi \left(jnx^{n-1} - kmx^{m-1}\right) = \xi (jnx^{n-m}-km)x^{m-1}.
$$
Clearly, for 
$$
x_1=\left(\frac{km}{jn}\right)^{1/(n-m)}
$$
we have $\phi'(x_1)=0$, and for any other $x>1$ we have $\phi'(x) \neq 0$. Furthermore, it is easily seen that for $x>1$ there is
only one possible value where $\phi''(x)=0$, namely the value
$$
x_2 = \left( \frac{km(m-1)}{jn(n-1)}\right)^{1/(n-m)}.
$$
Note that $x_2 < x_1$. This means that the interval $[\alpha, \beta]$ can be
partitioned into at most 3 subintervals, in all of which $\phi'(x)$ is monotonic, respectively. More precisely,
in the interval $(1, x_2]$ the function $\phi'$ is negative and
monotonic decreasing, in $[x_2,x_1]$ it is negative and monotonic
increasing, and in $[x_1, \infty)$ it is positive and monotonic increasing.\\

We have
$$
\frac{jnx_1^{n-1}}{kmx_1^{m-1}} = 1.
$$
Thus for any $x$ satisfying $x \geq x_1 \left(1+\eta\right)$ this implies
$$
\frac{jnx^{n-1}}{kmx^{m-1}} \geq \underbrace{\frac{jnx_1^{n-1}}{kmx_1^{m-1}}}_{=1}
\left(1+\eta\right)^{n-m} \geq 1 + \eta,
$$
and consequently
\begin{equation} \label{fi1}
\phi'(x)  \geq \left((1+\eta-1\right) \xi kmx^{m-1} \geq \xi \eta m \alpha^{m-1}.
\end{equation}
Similarly, for any $x$ satisfying $x \leq x_1 \left(1-\eta \right)$ we have
$$
\frac{jnx^{n-1}}{kmx^{m-1}} \leq \underbrace{\frac{jnx_1^{n-1}}{kmx_1^{m-1}}}_{=1}
\left(1-\eta\right)^{n-m} \leq 1 - \eta
$$
and consequently
\begin{equation} \label{fi2}
|\phi'(x)|  \geq \left|\left(1-\eta -1\right) \xi kmx^{m-1}\right| \geq \xi \eta m
\alpha^{m-1}. 
\end{equation}
Now we set 
$$
E = \left[x_1 \left(1-\eta \right), x_1 \left(1+\eta \right)\right]
$$
and
$$
I_1 = [A,x_2]\backslash E, \quad I_2 = [x_2, x_1] \backslash E, \qquad I_3 =
[x_1, B] \backslash E.
$$

Note that it is possible that some of these three intervals are empty, which is
no problem. Whenever $[\alpha, \beta]$ is completely contained in one of the intervals $I_1,
I_2$ or $I_3$ the derivative of $\phi(x)$ is monotonic in $[\alpha,\beta]$, and
by~\eqref{fi1} and~\eqref{fi2} for any $x\in [\alpha,\beta]$ we have
$$
|\phi'(x)| \geq \xi \eta m \alpha^{m-1}.
$$
Thus in this case by Lemma~\ref{vclemma}
$$
\left| \int_\alpha^\beta \cos (2 \pi \xi (j x^n - k x^m)) ~dx \right| \leq \frac{1}{\xi \eta m
\alpha^{m-1}}.
$$
Note also that
$$
\p \left([A,B] \backslash (I_1 \cup I_2 \cup I_3)\right) \leq \p(E) \leq 2 B \eta.
$$
This proves Lemma~\ref{lemma5}.\\

\section{Exponential inequality} \label{expi}

\begin{lemma} \label{lemma1}
Let $f$ be a function satisfying~\eqref{f}. Assume additionally that
\begin{equation} \label{fass}
\|f\| \geq \frac{N^{-1/4}}{\sqrt{2}}. 
\end{equation}
Then there exist numbers $c_A \geq 1$ and $N_0$ (depending only on $A$) such
that for any $M \geq 0,~N \geq N_0,$ and any $\delta > 0$ we have
\begin{eqnarray*}
& & \p \left( x \in [A,B]:~\left| \sum_{n=M+1}^{M+N} f(\xi x^{s_n}) \right| \geq
\delta c_A \|f\|^{1/4} \sqrt{N \log \log N} \right) \\
& \ll & \exp \left((1-\delta/2) \|f\|^{-1/2} \log \log N\right) + \delta^{-2}
N^{-16}.
\end{eqnarray*}
\end{lemma}

For the proof of Lemma~\ref{lemma1} we use a method of Takahashi~\cite{taka}, in
a refined form of Philipp~\cite{plt}. For simplicity of writing we assume that
$f$ is an even function, i.e. that it can be expanded into a pure cosine-series
$$
f(x) \sim \sum_{j=1}^\infty a_j \cos 2 \pi j x;
$$
the proof in the general case is exactly the same. Then by Lemma~\ref{ztlemma}
we have
\begin{equation} \label{aj}
|a_j| \leq \frac{1}{j}, \qquad j \geq 1.
\end{equation}
For any given $M \geq 0$, we write $(w_1, \dots, w_N)$ for the sequence $(s_{M+1}, \dots, s_{M+N})$. Set
$$
g(x) = \sum_{j=1}^{N^{38}} a_j \cos 2 \pi j x, \qquad r(x) = f(x)-g(x).
$$
Note that~\eqref{f} implies 
\begin{equation}  \label{*}
\|f\| \leq \|f\|_\infty \leq 1.
\end{equation}
Thus
\begin{equation} \label{ginf}
\|g\|_\infty \leq \|f\|_\infty + \var f \leq 3,
\end{equation}
by Lemma~\ref{ztlemma} and equations (1.25) and (3.5) of Chapter III
of~\cite{zt}.\\

Lemma~\ref{lemma1} will be deduced from Lemma~\ref{lemma1a} and
Lemma~\ref{lemma1b} below.

\begin{lemma} \label{lemma1a}
There exists a constant $\tilde{c}_A$ such that for any sufficiently large $N$
(depending only on $A$) and any $\tau>0$ satisfying 
\begin{equation} \label{lambda}
\tau \geq N^{-1/2} \qquad \textrm{and} \qquad 12 \tau  \left\lceil N^{1/8} \right\rceil \leq 1
\end{equation}
we have
\begin{eqnarray} \label{cs}
 \int_A^B \exp\left(\tau \sum_{n=1}^{N} g\left(\xi x^{w_n}\right) \right) ~dx \ll
e^{\tau^2 \tilde{c}_A \|f\| N}.
\end{eqnarray}
\end{lemma}

\begin{lemma} \label{lemma1b}
$$
\int_A^B \left(\sum_{n=1}^N r(\xi x^{w_n})\right)^2 ~dx \ll N^{-15}.
$$
\end{lemma}

\begin{corollary} \label{cor1}
Under the assumptions of Lemma~\ref{lemma1a}, we have
$$
\int_A^B \exp\left(\left| \tau \sum_{n=1}^{N} g(\xi x^{w_n}) \right|
\right) ~dx \ll e^{ \tau^2 \tilde{c}_A \|f\| N}.
$$ 
\end{corollary}
The corollary is obtained by using Lemma~\ref{lemma1a} also for the function
$-g(x)$ instead of $g(x)$.\\

\emph{Proof of Lemma~\ref{lemma1a}:}~We use the inequality
\begin{equation} \label{eque}
e^z \leq 1 + z + z^2, \qquad \textrm{for $|z| \leq 1$}.
\end{equation}
Set
\begin{equation} \label{H}
H = \left\lceil N^{1/8} \right\rceil,
\end{equation}
and 
$$
P = \max \{m \in \N:~Hm < N\}.
$$
For $1 \leq m < P$ set
\begin{equation} \label{hm}
U_m (x) = \sum_{n=Hm+1}^{H(m+1)} g(\xi x^{w_n}),
\end{equation}
and
$$
U_P(x) =  \sum_{n=HP+1}^{N} g(\xi x^{w_n}).
$$
Then
$$
\sum_{n=1}^N g(\xi x^{w_n}) = \sum_{m=1}^P U_m(x).
$$
By~\eqref{ginf} and the first inequality in~\eqref{lambda} we have
\begin{eqnarray} 
\int_A^B \exp \left( 4\tau U_{P}(x) \right) ~dx \leq e^1 \ll e^{\tilde{c}_A
\tau^2 N} \label{pest}
\end{eqnarray}
(the value of $\tilde{c}_A$ will be chosen later, but we can assume that $\tilde{c}_A \geq 1$). By the Cauchy--Schwarz inequality,~\eqref{cs} follows from~\eqref{pest}, together
with
\begin{equation} \label{u1}
\int_A^B  \exp \left( 4\tau \sum_{m:~1 \leq 2m < P} U_{2m}(x) \right) ~dx \ll
e^{\tilde{c}_A \tau^2 N},
\end{equation} 
and
\begin{equation} \label{u2}
\int_A^B \exp \left( 4\tau \sum_{m:~1 \leq 2m+1 < P} U_{2m+1}(x) \right) ~dx
\ll e^{\tilde{c}_A \tau^2 N}.
\end{equation}
The main idea of splitting the integral~\eqref{cs} into the parts~\eqref{u1} and~\eqref{u2} is that by separating the functions $U_m$ into two classes (those with even and those with odd index) there is also a separation of the corresponding values of $w_n$ in~\eqref{hm}. Consequently, two functions $U_{m_1}$ and $U_{m_2}$ which have both even or both odd index are ``almost'' independent, and $\int_A^B e^{U_{m_1}} e^{U_{m_2}} ~dx \approx \left(\int_A^B e^{U_{m_1}} ~dx\right)\left(\int_A^B e^{U_{m_2}} ~dx\right)$. Furthermore, the number of summands in the definition of $U_m$ is so small that we can use the approximation~\eqref{eque} for $\tau U_m$.\\

We will only prove~\eqref{u1}; the proof of~\eqref{u2} can be given in exactly
the same way. By~\eqref{ginf} and~\eqref{lambda} we have
$$
\left| 4\tau U_{2m}\right| \leq 12 \tau H \leq 1, \qquad 1 \leq 2m < P.
$$
By~\eqref{eque} this implies
\begin{eqnarray*}
\int_A^B  \exp \left( 4\tau \sum_{1 \leq 2m < P} U_{2m}(x) \right) ~dx & = &
\int_A^B  \prod_{1 \leq 2m < P} \exp \left( 4\tau U_{2m}(x) \right) ~dx \\
& \leq & \int_A^B \prod_{1 \leq 2m < P} \left( 1 + 4 \tau U_{2m}(x) + 16
\tau^2 U_{2m}(x)^2 \right) dx.
\end{eqnarray*}
For any $m,~1 \leq 2m < P$, using the standard trigonometric identity
\begin{equation} \label{standard}
\cos x \cos y = \frac{1}{2} \left( \cos(x+y) + \cos (x-y)\right),
\end{equation} 
we have
\begin{eqnarray}
16 \tau^2 U_{2m}^2 & = & 16 \tau^2 \left( \sum_{n=2Hm+1}^{H(2m+1)}
\sum_{j=1}^{N^{38}} a_j \cos \left(2 \pi \xi j x^{w_n}\right) \right)^2
\nonumber\\
& = & 16 \tau^2 \sum_{n_1,n_2=2Hm+1}^{H(2m+1)} ~\sum_{j_1,j_2=1}^{N^{38}}
\frac{a_{j_1} a_{j_2}}{2} \cos \left(2 \pi \xi (j_1 x^{w_{n_1}} + j_2
x^{w_{n_2}}) \right) \nonumber\\
& & + 16 \tau^2 \sum_{n_1,n_2=2Hm+1}^{H(2m+1)} ~\sum_{j_1,j_2=1}^{N^{38}}
\frac{a_{j_1} a_{j_2}}{2} \cos \left(2 \pi \xi (j_1 x^{w_{n_1}} - j_2
x^{w_{n_2}}) \right) \nonumber\\
& = & 16 \tau^2\sum_{n_1,n_2=2Hm+1}^{H(2m+1)} ~\sum_{j_1,j_2=1}^{N^{38}}
\frac{a_{j_1} a_{j_2}}{2} \cos \left(2 \pi \xi (j_1 x^{w_{n_1}} + j_2
x^{w_{n_2}}) \right) \label{Ux1}\\
& & + 16 \tau^2 \sum_{n_1,n_2=2Hm+1}^{H(2m+1)}
~\underbrace{\sum_{j_1,j_2=1}^{N^{38}}}_{*} \frac{a_{j_1} a_{j_2}}{2} \cos
\left(2 \pi \xi (j_1 x^{w_{n_1}} - j_2 x^{w_{n_2}}) \right) \label{Ux2} \\
& & + 16 \tau^2 \sum_{n_1,n_2=2Hm+1}^{H(2m+1)}
~\underbrace{\sum_{j_1,j_2=1}^{N^{38}}}_{**} \frac{a_{j_1} a_{j_2}}{2} \cos
\left(2 \pi \xi (j_1 x^{w_{n_1}} - j_2 x^{w_{n_2}}) \right). \label{Ux3}
\end{eqnarray}
Here the symbol ``$*$'' in~\eqref{Ux2} indicates that in this sum only those
values of $j_1,j_2$ are considered for which
$$
\max\{j_1,j_2\} \geq \min\{j_1,j_2\} A^{|n_1-n_2|},
$$
while the symbol ``$**$'' in~\eqref{Ux3} means that this sum is restricted to
those $j_1,j_2$ for which
\begin{equation} \label{ucond}
\max\{j_1,j_2\} < \min\{j_1,j_2\} A^{|n_1-n_2|}.
\end{equation}
We write $W_{2m}(x)$ for the sum of $4 \tau U_{2m}(x)$ plus the expressions
in~\eqref{Ux1} and~\eqref{Ux2}, and $V_{2m}$ for the expression in~\eqref{Ux3}.
Then 
\begin{equation} \label{wve}
1 + 4 \tau U_{2m}(x) + 16 \tau^2 U_{2m}(x)^2 = 1 + W_{2m}(x) + V_{2m}(x).
\end{equation}
For $V_{2m}(x)$ we have, using~\eqref{aj},~\eqref{ucond}, and the
Cauchy--Schwarz inequality,
\begin{eqnarray}
\|V_{2m}\|_\infty & \leq & 16 \tau^2 \sum_{n_1,n_2=2Hm+1}^{H(2m+1)}
\underbrace{\sum_{j_1,j_2=1}^{N^{38}}}_{\max\{j_1,j_2\} \geq \min\{j_1,j_2\}
A^{|n_1-n_2|}} \frac{a_{j_1} a_{j_2}}{2} \nonumber\\
& \leq &  64 \tau^2 \quad \sum_{2Hm+1 \leq n_1 \leq n_2 \leq H(2m+1)} \quad
\sum_{j_2=1}^{N^{38}} \quad \sum_{j_1 \geq j_2 A^{(n_2-n_1)}} \frac{a_{j_1}
a_{j_2}}{2} \nonumber\\
& \leq &  64\tau^2 H \quad \sum_{\ell = 0}^\infty \quad
\sum_{j_2=1}^{N^{38}} \quad \sum_{j_1 \geq j_2 A^{\ell}} \frac{a_{j_1}
a_{j_2}}{2} \nonumber\\
& \leq &  64\tau^2H \quad \left(\sum_{j=1}^{N^{38}}
\left(\frac{a_j}{\sqrt{2}}\right)^2 \right)^{1/2} \sum_{\ell = 0}^\infty \left(
\sum_{j \geq A^{\ell}} \frac{1}{j^2} \right)^{1/2} \nonumber\\
& \leq &  \tau^2H \|f\| \tilde{c}_A, \label{**}
\end{eqnarray}
where $\tilde{c}_A$ is a constant depending only on $A$. Thus by~\eqref{wve}
$$
1 + 4 \tau U_{2m}(x) + 16 \tau^2 U_{2m}(x)^2 \leq 1 + W_{2m}(x) + 
\tau^2 H \|f\| \tilde{c}_A.
$$
We can write the function $W_{2m}(x)$ as a sum of at most $3 H^2 N^{76}$
functions of the form
\begin{equation} \label{form}
\cos \left(2\pi \xi j x^\ell\right) \qquad \textrm{or} \qquad \cos \left(2\pi
\xi (j_1 x^{\ell_1} \pm j_2 x^{\ell_2})\right),
\end{equation}
all of which have coefficients bounded by $\max \{|a_j|,~j \geq 1\} \leq 1$. The
derivative of the arguments of the cosine-functions in~\eqref{form} is at most
\begin{equation} \label{der1}
4 \pi \xi N^{38} w_{H(2m+1)} x^{w_{H(2m+1)}-1} \qquad \textrm{for $x \in [A,B]$},
\end{equation}
and, on the other hand, by construction, this derivative is at least
\begin{equation} \label{der2}
2 \pi \xi w_{2Hm+1} x^{w_{2Hm+1}-1} \left( 1 - A^{-1} \right) \qquad \textrm{for $x \in
[A,B]$}
\end{equation}
(since functions having smaller derivative are collected in $V_{2m}$).
Furthermore, the second derivative is at least 
\begin{equation} \label{second}
2 \pi \xi w_{2Hm+1} \left(w_{2Hm+1}-1\right) x^{w_{2Hm+1}-2} \left( 1 - A^{-1} \right)
\qquad \textrm{for $x \in [A,B]$}.
\end{equation}
By~\eqref{wve} and~\eqref{**} we have
\begin{eqnarray*}
& & \int_A^B \prod_{1 \leq 2m < P} \left( 1 + 4 \tau U_{2m(x)} + 16 \tau^2
U_{2m}(x)^2 \right) dx \\
& \leq & \int_A^B \prod_{1 \leq 2m < P} \left( 1 + W_{2m}(x) +  \tau^2 H \|f\|
\tilde{c}_A \right) dx.
\end{eqnarray*}
For some $L$ let $i_1 \leq \dots \leq i_L$ be any numbers from the set $\{m:~1 \leq
2m < P\}$, and let $h_{i_1}(x), \dots, h_{i_L}(x)$ be functions of the form
\eqref{form} from $W_{2i_1}(x), \dots, W_{2i_L}(x)$, resp. Then by~\eqref{der1}
and~\eqref{der2} the product
\begin{equation} \label{prodf}
\prod_{\ell=1}^L h_{i_\ell}
\end{equation}
is a sum of cosine-functions with coefficients at most 1, such that the argument
of each cosine-function has derivative at least
\begin{eqnarray*}
& & 2 \pi \xi \left( w_{2Hi_L+1} x^{w_{2Hi_L+1}-1} \left(1-A^{-1}\right) -
\sum_{\ell=1}^{L-1} 2N^{38} w_{H(2i_\ell+1)} x^{w_{H(2i_\ell+1)}-1} \right) \\
& \geq & 2 \pi \xi w_{2Hi_L+1} \left(  x^{w_{2Hi_L+1}-1} \left(1-A^{-1}\right) - 2N^{38}
\sum_{\ell=1}^{L-1} x^{w_{H(2i_\ell+1)}-1} \right) \\
& \geq & 2 \pi \xi w_{2Hi_L+1} \left(  x^{w_{2Hi_L+1}-1} \left( \left(1-A^{-1}\right) - 2
N^{39} \underbrace{\frac{ x^{w_{H(2i_{L-1}+1)}-1}}{x^{w_{2Hi_L+1}-1}}}_{\leq x^{-H}} \right) \right) \\
& \geq & 2 \pi \xi w_{2Hi_L+1} \left(  x^{2Hi_L} \left(\left(1-A^{-1}\right) - 2N^{39}
x^{-H} \right) \right) \\
& \gg & A^{2Hi_L},
\end{eqnarray*}
for sufficiently large $N$, since by~\eqref{H}
$$
2N^{39} x^{-H} \leq \frac{(1-A^{-1})}{2}
$$
for sufficiently large $N$ (depending only on $A$). Similarly, using
\eqref{second}, it is seen that the second derivative of the argument of each
cosine-function which appears in~\eqref{prodf} is positive (for sufficiently large $N$). Thus by
Lemma~\ref{vclemma}
$$
\int_A^B \left( \prod_{\ell=1}^L h_{i_\ell} \right) ~dx \ll A^{-2Hi_L}.
$$
For any fixed $K$, there are in total at most 
$$
(3H^2N^{76})^{i_L}
$$
functions of the form~\eqref{prodf} for which $i_L=K$ (in other words, functions
which are composed from one function in $W_{2K}$ and at most one function from
$W_2, W_4, \dots, W_{2K-2}$). Furthermore, each of them has coefficient at most
1. Thus, using $1+x \leq e^x$ and
$$
\sum_{K=1}^{\lceil P/2 \rceil} \frac{(3H^2N^{76})^{2K}}{A^{2HK}} \leq \sum_{K=1}^\infty
\left(\frac{(3H^2N^{76})}{A^{H}} \right)^{2K} \ll 1 \qquad \textrm{for
sufficiently large $N$},
$$
we obtain
\begin{eqnarray*}
& &\int_A^B \prod_{1 \leq 2m < P} \left( 1 + W_{2m}(x) +  \tau^2 H \|f\| \tilde{c}_A
\right) dx \\
& \ll & \left(\prod_{1 \leq 2m < P} \left( 1 +  \tau^2 H \|f\| \tilde{c}_A \right)\right) \left( 1
+ \sum_{K=1}^{\lceil P/2\rceil} \frac{(3H^2N^{76})^{2K}}{A^{2HK}} \right) \\
& \ll & \prod_{1 \leq 2m < P} \exp \left(\tau^2 H \|f\| \tilde{c}_A \right) \\ 
& \ll & \exp \left(\tau^2 \tilde{c}_A \|f\| N \right).  
\end{eqnarray*}
This proves Lemma~\ref{lemma1a}.\\

\emph{Proof of Lemma~\ref{lemma1b}:}~ By Minkowski's inequality
\begin{equation} \label{erste}
\left( \int_A^B \left(\sum_{n=1}^N r(\xi x^{w_n})\right)^2 ~dx \right)^{1/2}
\leq \sum_{n=1}^N \left(\int_A^B \left(r(\xi x^{w_n})\right)^2 ~dx
\right)^{1/2}. 
\end{equation}
For $w_n$ fixed, using~\eqref{aj}, Lemma~\ref{lemma3} and the inequality of
arithmetic and geometric means we have
\begin{eqnarray*}
& & \int_A^B \left(r(\xi x^{w_n})\right)^2 ~dx \\
& = & \int_A^B \left( \sum_{j,k=N^{38}+1}^\infty \frac{a_j a_k}{2} \big( \cos (2 \pi
(j+k) \xi x^{w_n}) + \cos (2 \pi (j-k) \xi x^{w_n}) \big) \right) dx \\
& \leq & \sum_{j,k=N^{38}+1}^\infty \frac{1}{2jk} \left| \int_A^B \left( \cos (2 \pi
(j+k) \xi x^{w_n}) + \cos (2 \pi (j-k) \xi x^{w_n}) \right) dx \right| \\
& \leq & 2 \sum_{j=N^{38}+1}^\infty \sum_{\ell=0}^\infty
\frac{1}{2j(j+\ell)} \left| \int_A^B \left( \cos (2 \pi (2j+\ell) \xi x^{w_n}) + \cos (2 \pi
\ell \xi x^{w_n}) \right) dx \right| \\
& \leq & 2 \sum_{j=N^{38}+1}^\infty \left( \frac{1}{j^2} + \sum_{\ell=1}^\infty
\frac{1}{2j(j+\ell)} \left( \frac{1}{(2j+\ell)\xi w_n A^{w_n-1}} +
\frac{1}{\ell\xi w_n A^{w_n-1}}\right) \right) \\
& \leq & \frac{2}{N^{38}} + 2 \sum_{j=N^{38}+1}^\infty  \sum_{\ell=1}^\infty
\frac{1}{j(j+\ell) \ell} \frac{1}{\xi w_n A^{w_n-1}}  \\
& \leq & \frac{2}{N^{38}} + \frac{2}{\xi} \sum_{j=N^{38}+1}^\infty  \sum_{\ell=1}^\infty
\frac{1}{j(j+\ell)\ell} \\
& \leq & \frac{2}{N^{38}} + \frac{1}{\xi} \sum_{j=N^{38}+1}^\infty  \sum_{\ell=1}^\infty
\frac{1}{j^{3/2} \ell^{3/2}} \\
& \ll & N^{-17}.
\end{eqnarray*}
Together with~\eqref{erste} this proves Lemma~\ref{lemma1b}.\\

\emph{Proof of Lemma~\ref{lemma1}:}~ We use Corollary~\ref{cor1} for
$$
\tau = \tilde{c}_A^{-1/2} N^{-1/2} (\log \log N)^{1/2} \|f\|^{-3/4}
$$
(here $\tilde{c}_A$ is the constant from the statement of Corollary~\ref{cor1}). Then by~\eqref{fass} and~\eqref{*} we have
$$
\tau \geq  \tilde{c}_A^{-1/2} N^{-1/2} (\log \log N)^{1/2} \qquad \textrm{and} \qquad \tau \leq 2 \tilde{c}_A^{-1/2} N^{-5/16} (\log \log N)^{1/2},
$$
and condition~\eqref{lambda} is satisfied for sufficiently large $N$. Consequently
\begin{eqnarray*}
& & \int_A^B \exp\left(\left| \tilde{c}_A^{-1/2} N^{-1/2} (\log \log N)^{1/2}
\|f\|^{-3/4} \sum_{n=M+1}^{M+N} g(\xi x^{s_n}) \right| \right) ~dx \\
& \ll & 
e^{\|f\|^{-1/2} \log \log N},
\end{eqnarray*}
which implies that for arbitrary $\delta>0$ we have
\begin{equation} \label{le1a1}
\p \left(\left| \sum_{n=M+1}^{M+N} g(\xi x^{s_n}) \right| \geq (\delta/2)
\sqrt{\tilde{c}_A} \|f\|^{1/4} \sqrt{N \log \log N} \right) \ll e^{(1 -
\delta/2) \|f\|^{-1/2} \log \log N}.
\end{equation}
By Lemma~\ref{lemma1b}, Markov's inequality and~\eqref{fass} we have
\begin{equation} \label{le1a2}
\p \left(\left| \sum_{n=M+1}^{M+N} r(\xi x^{s_n}) \right| \geq (\delta/2)
\sqrt{\tilde{c}_A} \|f\|^{1/4} \sqrt{N \log \log N} \right) \ll \|f\|^{-1/2}
\delta^{-2} N^{-17} \leq \delta^{-2} N^{-16}.
\end{equation}
Combining~\eqref{le1a1} and~\eqref{le1a2} we finally obtain
$$
\p \left(\left| \sum_{n=M+1}^{M+N} f(\xi x^{s_n}) \right| \geq \delta
\underbrace{\sqrt{\tilde{c}_A}}_{=:c_A} \|f\|^{1/4} \sqrt{N \log \log N}
\right) \ll e^{(1 - \delta/2) \|f\|^{-1/2} \log \log N} + \delta^{-2} N^{-16},
$$
which proves Lemma~\ref{lemma1}.

\section{Maximal inequality} \label{maxi}

\begin{lemma} \label{lemmamax}
For any sufficiently large $m$ (depending only on $A$) we have the following:
Let $f$ be a function satisfying~\eqref{f} and
\begin{equation} \label{a*}
\|f\| \geq \frac{2^{-m/4}}{\sqrt{2}}.
\end{equation}
Then for the number $c_A$ from Lemma~\ref{lemma1} and any $\gamma \geq 1$ we have
\begin{eqnarray*}
& & \p \left( x \in [A,B]:~\max_{1 \leq M \leq 2^m} \left| \sum_{n=1}^{M} f(\xi
x^{s_n}) \right| \geq 84 \gamma c_A \|f\|^{1/4} \sqrt{2^m \log \log 2^m} \right)
\\
& \ll & \exp \left(-\gamma \|f\|^{-1/2} \log \log 2^m \right) + 2^{-3m}.
\end{eqnarray*}
\end{lemma}

\begin{corollary} \label{co3}
For any sufficiently large $N$ (depending only on $A$) we have the following:
Let $f$ be a function satisfying~\eqref{f}. Assume additionally that
\eqref{fass} holds. Then for the number $c_A$ from Lemma~\ref{lemma1} and any $\gamma \geq 1$ we have
\begin{eqnarray*}
& & \p \left( x \in [A,B]:~\max_{1 \leq M \leq N} \left| \sum_{n=1}^{M} f(\xi
x^{s_n}) \right| \geq 119 \gamma c_A \|f\|^{1/4} \sqrt{N \log \log N} \right) \\
& \ll & \exp \left(-\gamma \|f\|^{-1/2} \log \log N\right) + N^{-3}.
\end{eqnarray*}
\end{corollary}

\emph{Proof of Lemma~\ref{lemmamax}:~}For the proof of Lemma~\ref{lemmamax} we use
a classical dyadic decomposition method, which is frequently used for proving
maximal inequalities in probability theory and probabilistic number theory (see,
for example,~\cite{baker,gako}). By Lemma~\ref{lemma1} for the complete sum
$\sum_{n=1}^{2^m}f(\xi x^{s_n})$ we have
\begin{eqnarray}
& & \p \left(\left| \sum_{n=1}^{2^m} f(\xi x^{s_n}) \right| \geq 84 \gamma c_A
\|f\|^{1/4} \sqrt{N \log \log N} \right) \nonumber\\
& \ll & \exp \left((1-84 \gamma/2) \|f\|^{-1/2} \log \log 2^m\right) + 2^{-16m} \nonumber\\
& \ll & \exp \left(-41\gamma \|f\|^{-1/2} \log \log 2^m\right) + 2^{-16m}. 
\label{totalsum}
\end{eqnarray}
Any number $M < 2^m$ can be written in dyadic representation
$$
M = \ve_0 + 2 \ve_1 + 4 \ve_2 + \dots + 2^{m-1} \ve_{m-1} \qquad \textrm{for
digits $\ve_0, \ve_1, \dots, \ve_{m-1}$}.
$$
Writing $\mathcal{S}$ for the set of those numbers $M,~1 \leq M \leq 2^m-1$, for which
$\ve_0 = 0, \ve_1 = 0, \dots, \ve_{m/4}=0$, then by~\eqref{a*}, for sufficiently large $m$,
\begin{eqnarray}
& & \p \left(\max_{1 \leq M \leq 2^m-1} \left| \sum_{n=1}^{M} f(\xi x^{s_n})
\right| \geq 84 \gamma c_A \|f\|^{1/4} \sqrt{2^m \log \log 2^m} \right)
\nonumber\\
& \leq & \p \left(\max_{M \in \mathcal{S}} \left| \sum_{n=1}^{M} f(\xi x^{s_n}) \right|
\geq 83 \gamma c_A \|f\|^{1/4} \sqrt{2^m \log \log 2^m} \right). \label{uksum}
\end{eqnarray}
For a set $U(K)$ containing $2^K$ consecutive elements of $\{1, \dots,
2^{m}-1\}$ for some $K,~ m/4 \leq K \leq m-1$, we have, using the fact that
$$
\log \log 2^K \geq \frac{\log \log 2^m}{2} \qquad \textrm{and} \qquad 1 - 5(m-K)
\geq 4(m-K)
$$
for sufficiently large $m$, and using Lemma~\ref{lemma1} for $\delta=10 (m-K) \gamma$,
\begin{eqnarray*}
& & \p \left(\left| \sum_{n \in U(K)} f(\xi x^{s_n}) \right| \geq 10
(m-K)2^{(K-m)/2} \gamma c_A \|f\|^{1/4} \sqrt{2^m \log \log 2^m} \right) \\
& \leq & \p \left(\left| \sum_{n \in U(K)} f(\xi x^{s_n}) \right| \geq 10 (m-K)
\gamma c_A \|f\|^{1/4} \sqrt{2^K \log \log 2^K} \right) \\
& \ll & \exp \left((1- 10(m-K) \gamma/2) \|f\|^{-1/2} \log \log 2^K \right) +
2^{-16K} \\
& \ll & \exp \left(-2(m-K) \gamma \|f\|^{-1/2} \log \log 2^m \right) +
2^{-16m/4} \\
& \ll & 2^{-2(m-K)} \exp \left(-\gamma \|f\|^{-1/2} \log \log 2^m \right) +
2^{-4m},
\end{eqnarray*}
provided $m$ is sufficiently large. To be able to represent every set $\{1, \dots,
M\}$ for $M \in \mathcal{S}$ as a disjoint union of at most one set of cardinality $2^K$ for each $K \in \{K:~m/4 \leq K \leq m-1\}$, we need in total $2^{m-K}$ sets of cardinality $2^K$, for each $m/4 \leq K
\leq m-1$. Thus, using
$$
\sum_{k=1}^\infty 10 k 2^{-k/2} \leq 83,
$$
we have
\begin{eqnarray*}
\left| \sum_{n=1}^M f(\xi x^{s_n}) \right| & < & \sum_{K=m/4}^{m-1} 10
(m-K)2^{(K-m)/2} \gamma c_A \|f\|^{1/4} \sqrt{2^m \log \log 2^m} \\
& < & 83 \gamma c_A \|f\|^{1/4} \sqrt{2^m \log \log 2^m}
\end{eqnarray*}
for all $M \in \mathcal{S}$, except for a set $x \in [A,B]$ of measure at most
\begin{eqnarray*}
& \ll & \sum_{K=m/4}^{m-1} 2^{m-K} \left( 2^{-2(m-K)} \exp \left(-\gamma
\|f\|^{-1/2} \log \log 2^m \right) + 2^{-4m} \right) \\
& \ll & \exp \left(-\gamma \|f\|^{-1/2} \log \log 2^m \right) + 2^{-3m},
\end{eqnarray*}
provided $m$ is sufficiently large. Together with~\eqref{totalsum} and
\eqref{uksum} this proves Lemma~\ref{lemmamax}.\\

\emph{Proof of Corollary~\ref{co3}:~} Write $\hat{N}$ for the smallest number
$\geq N$, which is a power of $2$. Then $2N > \hat{N} \geq N$. Using
Lemma~\ref{lemmamax} we have
\begin{eqnarray*}
\max_{1 \leq M \leq N} \left| \sum_{n=1}^{M} f(\xi x^{s_n}) \right| & \leq & \max_{1 \leq M \leq \hat{N}} \left| \sum_{n=1}^{M} f(\xi x^{s_n}) \right| \\
& < & 84 \gamma c_A \|f\|^{1/4} \sqrt{\hat{N} \log \log \hat{N}}  \\
& \leq & 119 \gamma c_A \|f\|^{1/4} \sqrt{N \log \log N} 
\end{eqnarray*}
for sufficiently large $N$, except for a set $x \in [A,B]$ of probability at
most
$$
\exp \left(-\gamma \|f\|^{-1/2} \log \log \hat{N} \right) +
\left(\hat{N}\right)^{-3} \ll \exp \left(-\gamma \|f\|^{-1/2} \log \log N
\right) + N^{-3}.
$$

\section{The law of the iterated logarithm for functions having small
$L^2$-norm} \label{smallr}

\begin{lemma} \label{lemmar}
Let $f$ be a function of bounded variation satisfying~\eqref{f}. For some
$d \geq 1$, let $p$ denote the $d$-th partial sum of the Fourier series of $f$,
and let $r$ denote the remainder term $f-p$. Then for the constant $c_A$ from
Lemma~\ref{lemma1} we have
$$
\limsup_{N \to \infty} \frac{\left| \sum_{n=1}^N r(x^{s_n}) \right|}{\sqrt{N
\log \log N}} \leq 238 c_A d^{-1/8} \quad \textup{a.e.}
$$
\end{lemma}

\emph{Proof of Lemma~\ref{lemmar}:~} Using Lemma~\ref{ztlemma} we have
$$
\|r\| \leq \sum_{j=d+1}^\infty \frac{2}{j^2} \leq \frac{2}{\sqrt{d}}, 
$$
and, by~\eqref{f}, we also have $\|r\| \leq 1$. Setting
$$
D_m := \left( x \in [A,B]:~\max_{1 \leq M \leq 2^m} \left| \sum_{n=1}^M r(\xi
x^{s_n}) \right| > 84 (2d^{-1/8}) c_A \sqrt{2^m \log \log 2^m} \right)
$$
and using Lemma~\ref{lemmamax} for $\gamma=2d^{-1/8}\|r\|^{-1/4}\geq 1.6$ we
obtain
\begin{eqnarray*}
\p(D_m) & \ll & \exp \left(-1.6 \log \log 2^m \right) + 2^{-3m} \\
& \ll & \frac{1}{m^{1.6}}.
\end{eqnarray*}
Thus 
$$
\sum_{m=1}^\infty \p(D_m) < \infty,
$$
and by the Borel-Cantelli lemma with probability 1 only finitely many
events $D_m$ occur. Thus since $238 > \sqrt{2} \cdot 2 \cdot 84$ there are, also
with probability 1, only finitely many $N$ for which
$$
\left| \sum_{n=1}^N r(\xi x^{s_n}) \right| > 238 c_A d^{-1/8} \sqrt{N \log \log N},
$$
which proves the lemma.

\section{The law of the iterated logarithm for the discrepancy for small intervals}
\label{smalln}

In the present section we will prove a bounded law of the iterated logarithm for
a modified version of the discrepancy, which only takes into account ``small'' intervals. More precisely, for an integer $R
\geq 1$ and a sequence $(z_1, \dots, z_N) \in [0,1)^N$ set
\begin{equation} \label{dndef}
D_N^{(\leq 2^{-R})} (z_1, \dots, z_N) := \sup_{a \in \Z, 0 \leq a < 2^R,~0
\leq b \leq 2^{-R}} \quad \left| \frac{1}{N} \sum_{n=1}^N \mathbf{I}_{[a 2^{-R}
,a2^{-R} + b)} (z_n) \right|
\end{equation}
(the functions $\mathbf{I}$ were defined in~\eqref{I}). In other words, the discrepancy $D_N^{(\leq 2^{-R})}$ considers only ``small''
intervals (those of length $\leq 2^{-R}$), which have their left corner in a
point of the form $a2^{-R}$ for some $a \in \{0, \dots, 2^R-1\}$. Furthermore, we set 
$$
D_N^{(\geq 2^{-R})} (z_1, \dots, z_N) := \max_{a,b \in \Z, 0 \leq a < b \leq
2^R} \quad \left| \frac{1}{N} \sum_{n=1}^N \mathbf{I}_{[a2^{-R},b2^{-R})} (z_n)
\right|
$$
and
$$
{D_N^*}^{(\geq 2^{-R})} (z_1, \dots, z_N) := \max_{a \in \Z, 0 < a \leq 2^R}
\quad \left| \frac{1}{N} \sum_{n=1}^N \mathbf{I}_{[0,a2^{-R})} (z_n) \right|.
$$
It is easily seen that always
\begin{equation} \label{discrepancies}
{D_N^*}^{(\geq 2^{-R})} \leq D_N^* \leq D_N \leq D_N^{(\geq 2^{-R})} + 3
D_N^{(\leq 2^{-R})}.
\end{equation}
The idea to split the discrepancies $D_N^*$ and $D_N$ in this way to obtain
precise metric discrepancy results is due to Fukuyama~\cite{ft}.

\begin{lemma} \label{lemmasmint}
For any positive integer $R$ we have for almost all $x \in [A,B]$
$$
\limsup_{N \to \infty} \frac{\sqrt{N} D_N^{(\leq 2^{-R})}(\{\xi
x^{s_n}\})}{\sqrt{\log \log N}} \leq 10^7 c_A R^{-1},
$$
where $c_A$ is the constant from Lemma~\ref{lemma1}.
\end{lemma}

We use a dyadic decomposition of the unit interval, which was also used in
\cite{plt}. For simplicity we will only consider the case $a=0$, i.e.
$$
\limsup_{N \to \infty} \frac{\sup_{0 \leq b \leq 2^{-R}} \left| \sum_{n=1}^N
\mathbf{I}_{[0,b)} (\xi x^{s_n}) \right|}{\sqrt{N \log \log N}} \leq 10^7 c_A
R^{-1}
$$
for almost all $x \in [A,B]$. The proof for the other possible values of $a$, that is for
$1 \leq a \leq 2^R$, can be given in exactly the same way. This means that the
exceptional set in Lemma~\ref{lemmasmint} is a finite union of sets of measure
zero, and consequently also has zero measure.\\

For $N \geq 1$ we set 
$$
E_N = \left(\sup_{0 \leq b \leq 2^{-R}} \left| \sum_{n=1}^N \mathbf{I}_{[0,b)}
(\xi x^{s_n}) \right| \geq 10^7 c_A R^{-1} \sqrt{N \log \log N} \right),
$$
and for $m \geq 1$ we set
$$
F_m = \left(\max_{1 \leq N \leq 2^{2m}} \quad \sup_{b \in \Z, 1 \leq b \leq
2^{m-R}} \left| \sum_{n=1}^N \mathbf{I}_{[0,b2^{-m})} (\xi x^{s_n}) \right| \geq
10^6 c_A R^{-1} \sqrt{2^{2m} \log \log 2^{2m}} \right).
$$
Every interval $[0,b), ~0 \leq b \leq 2^{-R},$ can be written as the union of an
interval of the form $[0,j2^{-m})$ for some appropriate $j \in \Z,~1 \leq j < 2^{m-R},$ and an
interval $B$ of length at most $2^{-m}$. For any $x$ from the complement of
$F_m$ we have for any $N,~1 \leq N \leq 2^{2m}$, and any such interval $B$, 
\begin{eqnarray*}
\left| \sum_{n=1}^N \mathbf{I}_{B} (\xi x^{s_n}) \right| & < & 2 \cdot 10^6 c_A R^{-1}
\sqrt{2^{2m} \log \log 2^{2m}} + N 2^{-m} \\
& < & \left(2 \cdot 10^6+1\right) c_A R^{-1} \sqrt{2^{2m} \log \log 2^{2m}},
\end{eqnarray*}
provided $m$ is sufficiently large (depending on $A$ and $R$). Consequently for
any sufficiently large $N$ satisfying $2^{2m-2} \leq N \leq 2^{2m}$ for some $m$
we have, for any $x$ from the complement of $F_m$,
\begin{eqnarray*}
\sup_{0 \leq b \leq 2^{-R}} \left| \sum_{n=1}^N \mathbf{I}_{[0,b)} (\xi x^{s_n})
\right| & < & (3 \cdot 10^6+1) c_A R^{-1} \sqrt{2^{2m} \log \log 2^{2m}} \\
& \leq & 10^7 c_A R^{-1} \sqrt{N \log \log N}.
\end{eqnarray*}
Thus for sufficiently large $m$ we have
$$
\bigcup_{2^{2m-2} \leq N \leq 2^{2m}} E_N \subset F_m,
$$
and hence
\begin{equation} \label{fconv}
\sum_{m=1}^\infty F_m < \infty \qquad \textrm{implies} \qquad \sum_{N=1}^\infty
E_N < \infty.
\end{equation}
Writing $b$ in binary expansion, it is easily seen that for any possible number
$1 \leq b < 2^{m-R}$ the interval $[0,b2^{-m})$ can be written as the disjoint
union of at most one interval of length $2^{-R-1}$, at most one interval of
length $2^{-R-2}$, etc., and at most one interval of length $2^{-m}$. Furthermore,
to be able to represent all possible intervals $[0,b2^{-m})$ we need exactly
\begin{equation} \label{anzahl}
2^{k-R} \quad \textrm{intervals of length} \quad 2^{-k}, \qquad \textrm{for any~}k \in \{R+1, \dots, m\}.
\end{equation}
Let $f$ be the indicator function of an interval of length $2^{-k}$ for some $k
\leq m$. Then 
$$
\frac{2^{-k/2}}{\sqrt{2}} \leq \|f\| \leq 2^{-k/2}
$$
and
$$
\|f\| \geq \frac{2^{-2m/4}}{\sqrt{2}},
$$
Consequently, using Lemma~\ref{lemmamax} with $\gamma = 2 k$ we obtain 
\begin{eqnarray*}
& & \p \left( \max_{1 \leq N \leq 2^{2m}} \left|\sum_{n=1}^N f(\xi x^{s_n})
\right| \geq 84 (2 k) 2^{-k/8} c_A \sqrt{2^{2m} \log \log 2^{2m}} \right) \\
& \ll & \exp \left( -2k \log \log 2^{2m} \right) + 2^{-6m} \\
& \ll & \left(\frac{1}{m}\right)^{2k} + 2^{-6m}.
\end{eqnarray*}
It can be shown that
$$
\sum_{k=R+1}^m  84 (2 k) 2^{-k/8} \leq 250000 R^{-1}, \qquad \textrm{for any $R \geq 1$}.
$$
Thus, using~\eqref{anzahl}, we see that for any $b,~1 \leq b < 2^{m-R},$ we have
\begin{eqnarray*}
\max_{1 \leq N \leq 2^{2m}} \left| \sum_{n=1}^N \mathbf{I}_{[0,b2^{-m})} (\xi
x^{s_n}) \right| & \leq & \sum_{k=R+1}^m  84 (2 k) 2^{-k/8} c_A \sqrt{2^{2m}
\log \log 2^{2m}} \\
& \leq & 250000 c_A R^{-1} \sqrt{2^{2m}
\log \log 2^{2m}},
\end{eqnarray*}
except for a set of measure at most
\begin{eqnarray*}
& \ll & \sum_{k=R+1}^m 2^{k-R} \left(\left(\frac{1}{m}\right)^{2k} +
2^{-6m}\right) \\
& \ll & m^{-2}.
\end{eqnarray*}
Furthermore we also have for the full interval $[0,2^{-R})$, again by Lemma~\ref{lemmamax},
\begin{eqnarray*}
\p \left( \max_{1 \leq N \leq 2^{2m}} \left|\sum_{n=1}^N \mathbf{I}_{[0,2^{-R})} (\xi x^{s_n})
\right| \geq 250000 c_A R^{-1}\sqrt{2^{2m} \log \log 2^{2m}} \right) & \ll & m^{-2}.
\end{eqnarray*}
Thus
$$
\p (F_m) \ll m^{-2},
$$
which by~\eqref{fconv} and the Borel-Cantelli lemma implies that with
probability 1 only finitely many events $E_N$ occur. This proves Lemma
\ref{lemmasmint}.

\section{Martingale approximation} \label{sectmar}

Throughout this section we assume that the number $d$ is fixed; also throughout this
section the constants $c$ and the implied constants in ``$\ll$'' and ``$\gg$'' may depend on $d$. We will exclude the trivial case $\|p\|=0$, which is equivalent to $p \equiv 0$.

\begin{lemma} \label{lemma2}
Let $p(x)$ be a trigonometric polynomial. Then for all numbers $N$ which can be
written in the form
\begin{equation} \label{Nform}
N = \sum_{i=1}^M (i^4+i) \qquad \textrm{for some $M$}
\end{equation}
we have
$$
\sup_{t \in \R} \left| \p \left( x \in [A,B]:~\frac{\sum_{n=1}^N p(\xi x^{s_n})}{\|p\| \sqrt{N}} <
t \right) - \Phi(t) \right| \ll \frac{\log N}{N^{1/25}},
$$
where $\Phi$ denotes the standard normal distribution function
\end{lemma}

\begin{lemma} \label{lemma2lil}
Let $p(x)$ be a trigonometric polynomial, and let $(N_k)_{k \geq 1}$ be the sequence of numbers which can be
written in the form~\eqref{Nform}. Then
$$
\limsup_{k \to \infty} \frac{\left| \sum_{n=1}^{N_k} p(\xi x^{s_n}) \right|}{\sqrt{N_k \log \log N_k}} = \sqrt{2} \|p\| \qquad \textup{a.e.}
$$
\end{lemma}

The crucial ingredient in the proofs of Lemma~\ref{lemma2} and Lemma~\ref{lemma2lil}, which will be given simultaneously, are the following results of Strassen and of Heyde and Brown~\cite{hbo}, which are a consequence of a martingale version of the
Skorokhod representation theorem due to Strassen~\cite{str}. For Lemma~\ref{lstr} (which is used to prove Theorem~\ref{th1}) we use the formulation from~\cite[Lemma~2.1]{a2}, for Lemma~\ref{heyde} (which is used to prove Theorem~\ref{th2}) we use the formulation from~\cite[Theorem~A]{be78}.

\begin{lemma} \label{lstr}
Let $Y_1, Y_2, \dots$ be a
martingale difference sequence with finite fourth moments, let
$V_M = \sum_{i=1}^M \mathbb{E} (Y_i^2 | Y_1, \dots, Y_{i-1})$ and assume that $V_1 >0$ and $V_M \to \infty$. Let $(b_M)_{M \geq 1}$ be any sequence of positive numbers such that
$$
\lim_{M \to \infty} \frac{V_M}{b_M} = 1 \qquad \textup{a.s.},
$$
and
$$
\sum_{M=1}^\infty \frac{(\log b_M)^{10}}{b_M^2} \E Y_M^4 < \infty.
$$
Then
$$
\limsup_{M \to \infty} \frac{\sum_{i=1}^M Y_i}{\sqrt{b_M \log \log b_M}} = \sqrt{2} \qquad \textup{a.s.}
$$
\end{lemma}

\begin{lemma} \label{heyde}
Let $Y_1, Y_2, \dots$ be a
martingale difference sequence with finite fourth moments, let
$V_M = \sum_{i=1}^M \mathbb{E} (Y_i^2 | Y_1, \dots, Y_{i-1})$ and
let $(b_M)_{M \geq 1}$ be any sequence of positive numbers. Then
\begin{equation*}
\sup_{t \in \R} \left| \p \left( \frac{Y_1 + \dots + Y_M}{\sqrt{b_M}} <
t \right) - \Phi(t) \right| \leq K \left( \frac{\sum_{i=1}^M
\mathbb{E}Y_i^4 + \mathbb{E} \left( (V_M - b_M)^2 \right)}{b_M^2}
\right)^{1/5},
\end{equation*}
where $K$ is an absolute constant.
\end{lemma}

\emph{Proof of Lemma~\ref{lemma2}:}~ We use an argument based on approximation
by martingale differences, which was already used in~\cite{a2,a1}. This method
was originally developed by Berkes~\cite{bew,beo2,beo1} and Philipp and
Stout~\cite{ps}. Since in our case the functions which we want to approximate are not periodic, we have to
construct an increasing sequence of space-inhomogeneous discrete sigma-algebras for the
approximation.\\

For simplicity of writing we assume that $p$ is an even function; the proof in
the general case is exactly the same. Then we can write $p$ in the form
$$
p(x) = \sum_{j=1}^d a_j \cos 2\pi j x.
$$
For simplicity of writing we will also assume that $\|p\|\leq 1$ and $|a_j| \leq 1,~j \geq 1$.\\

We subdivide the set of positive integers consecutively into blocks $\Delta_i$
(``large blocks'') and $\Delta_i'$ (``small blocks''), in such a way that 
\begin{itemize}
 \item the block $\Delta_i$ contains $i^4$ elements, for $i \geq 1$.
 \item the block $\Delta_i'$ contains $i$ elements, for $i \geq 1$.
 \item elements of $\Delta_i$ are smaller than elements of $\Delta_i'$, for $i \geq
1$.
 \item elements of $\Delta_i'$ are smaller than elements of $\Delta_{i+1}$, for $i
\geq 1$.
 \item $\bigcup_{i\geq1} (\Delta_i \cup \Delta_i') = \N.$
\end{itemize}
Assume that $N$ is of the form~\eqref{Nform} for some $M$. Then by construction
$$
\{1, \dots, N\} = \bigcup_{i=1}^M \left( \Delta_i \cup \Delta_i'\right).$$
We write $min(i)$ and $max(i)$ for the smallest resp. largest
element of $\Delta_i$, and set
$$
m(i) = \left\lceil \log_2 \left(i^6 w_{max(i)} A^{w_{max(i)}}\right)
\right\rceil, \qquad  1 \leq i \leq M.
$$
Note that 
\begin{equation} \label{igl}
w_{\min(i)} - w_{\max(i-1)} \geq \min(i) - \max(i-1) \geq i-1, \qquad 1 \leq i \leq M.
\end{equation}
We write $G_i$ for the set of intervals of the form
$$
H_j^{(i)} = \left[A + j 2^{-m(i)}, A + (j+1) 2^{-m(i)} \right), \quad j \in \{0,
\dots, 2^{m(i)} - 1\}, \quad 1 \leq i \leq M.
$$
In other words, $G_i$ is a partition of $[A,B]$ into $2^{m(i)}$ subintervals of equal length. Write $x_j^{(i)}$ for the smallest number in $H_j^{(i)}$, that is
$$
x_j^{(i)} = A + j 2^{-m(i)},
$$
and split every interval $H_j^{(i)}$ into 
$$
2^{\left\lceil \log_2 \left(w_{\max(i)} \left(x_j^{(i)}/A\right)\right)
\right\rceil}
$$
pieces of equal length. Let $\F_i$ denote the sigma-algebra generated by all
these sets. Then $(\F_i)_{1 \leq i \leq M}$ is an increasing family of sigma-algebras.
For any $i$, any number $x \in [A,B]$ is contained in an atom of $\F_i$ which
has length between
\begin{eqnarray} \label{lengex}
\frac{2^{-m(i)}}{2} \left(\frac{A}{x}\right)^{w_{\max(i)}} \qquad \textrm{and}
\qquad 2^{-m(i)}\left(\frac{A}{x - 2^{-m(i)}}\right)^{w_{\max(i)}}.
\end{eqnarray}
For $1 \leq i \leq M$ we set
$$
T_i = \sum_{n \in \Delta_i} p(\xi x^{w_n}), \qquad T_i' = \sum_{n \in \Delta_i'}
p(\xi x^{w_n})
$$
and
$$
Y_i = \E \left( T_i \Big| \F_i \right) - \E \left( T_i \Big| \F_{i-1} \right).
$$
Then
$$
\frac{1}{\sqrt{N}} \sum_{n=1}^N p(\xi x^{s_n}) = \frac{1}{N} \sum_{i=1}^M \left( T_i + T_i'\right).
$$
Furthermore, $Y_i$ is a discrete function, which by construction is constant on the atoms of $\F_i$, and
$$
\E \left( Y_i \big| \F_{i-1} \right) = 0.
$$

In other words, $(Y_i)_{1 \leq i \leq M}$ is a martingale difference. Let
$[\alpha,\beta)$ be any atom of $\F_{i-1}$. Then by~\eqref{lengex}
\begin{eqnarray} 
\beta - \alpha \geq \frac{2^{-m(i-1)}}{2}
\left(\frac{A}{\alpha}\right)^{w_{\max(i-1)}} & \geq & \frac{1}{2 (i-1)^6
w_{\max(i-1)} A^{w_{\max(i-1)}}}
\left(\frac{A}{\alpha}\right)^{w_{\max(i-1)}} \nonumber\\
& \gg & \frac{1}{i^6 w_{\max(i-1)} \alpha^{w_{\max(i-1)}}}.
\label{lengex2}
\end{eqnarray}
Thus by Lemma~\ref{lemma3} and~\eqref{igl}
\begin{eqnarray*}
\left|\frac{1}{\beta-\alpha} \int_\alpha^\beta T_i(x) ~dx\right| & \leq &
\frac{1}{\beta-\alpha} ~\sum_{n \in \Delta_i} ~\sum_{j=1}^d |a_j| \left| \int_\alpha^\beta
 \cos 2 \pi j \xi x^{w_n} ~dx \right| \\
& \ll & \frac{i^6 w_{\max(i-1)} \alpha^{w_{\max(i-1)}}}{w_{\min(i)}
\alpha^{w_{\min(i)}}} \\
& \ll & i^6 A^{-i},
\end{eqnarray*}
and consequently
\begin{equation} \label{T_iab}
\E \left(T_i \big| \F_{i-1} \right) \ll i^6 A^{-i}.
\end{equation}
Now, let $[\alpha,\beta)$ denote an atom of $\F_i$. By~\eqref{lengex} we have
\begin{equation} \label{lengex3}
\beta - \alpha \leq 2^{-m(i)} \left( \frac{A}{\alpha-2^{-m(i)}}\right)^{w_{\max(i)}} \ll \frac{1}{i^6 w_{\max(i)} (\alpha-2^{-m(i)})^{w_{\max(i)}}}.
\end{equation}
The derivative of $p(\xi x^{w_n})$ on $[\alpha,\beta]$ is bounded by
$$
|p'(\xi x^{w_n})| \leq \sum_{j=1}^d |2 \pi j \xi w_n \beta^{w_n-1}| \ll w_n \left(\alpha+2^{-m(i)}\right)^{w_n}.
$$
By the definition of $m(i)$ it is easily seen that
$$
\left(\frac{\alpha+2^{-m(i)}}{\alpha-2^{-m(i)}}\right)^{w_{\max(i)}} \ll 1.
$$
Thus by~\eqref{lengex3} for any $n \in \Delta_i$ the fluctuation of $p(\xi x^{w_n})$ on $[\alpha,\beta)$ is bounded by $\ll i^{-6}$. Therefore, together with~\eqref{T_iab}, we obtain
\begin{eqnarray}
\left| T_i - Y_i \right| \ll |\Delta_i| i^{-6} \ll i^{-2} \label{yzest}
\end{eqnarray}
(here, and in the sequel, we write $|\cdot|$ for the number of elements of a set).\\

Next we have to calculate the conditional variances $\E (Y_i^2 | \F_{i-1})$. By
\eqref{yzest}
\begin{eqnarray}
\left| \E (Y_i^2 | \F_{i-1}) - \E (T_i^2 | \F_{i-1}) \right| & \leq & \E
\left(|(Y_i + T_i)(Y_i - T_i)| \Big| \F_{i-1}\right) \nonumber\\
& \ll & \|Y_i + T_i\|_\infty \|Y_i - T_i\|_\infty \nonumber\\
& \ll & |\Delta_i| i^{-2} \ll i^2, \label{tiyi}
\end{eqnarray}
and thus we can reduce the problem to estimating $\E (T_i^2 | \F_{i-1})$. Using
\eqref{standard}, we have
\begin{eqnarray*}
T_i^2 & = & \left( \sum_{n \in \Delta_i} \sum_{j=1}^d a_j \cos (2 \pi j \xi x^{w_n})
\right)^2 \\
& = & \frac{1}{2} \sum_{n_1,n_2 \in \Delta_i} ~\sum_{j_1,j_2=1}^d a_{j_1} a_{j_2}
\left( \cos (2 \pi \xi (j_1 x^{w_{n_1}} + j_2 x^{w_{n_2}})) +  \cos (2 \pi \xi (j_1
x^{w_{n_1}} - j_2 x^{w_{n_2}})) \right).
\end{eqnarray*}
In the above sum, for $n_1=n_2$ and $j_1=j_2$ we have $\cos (2 \pi \xi (j_1
x^{w_{n_1}} - j_2 x^{w_{n_2}}))=1$, and thus
$$
\frac{1}{2} \underbrace{\sum_{n_1,n_2 \in \Delta_i} ~\sum_{j_1,j_2=1}^d}_{(j_1,n_1)
= (j_2,n_2)} a_{j_1} a_{j_2} \cos (2 \pi \xi (j_1 x^{w_{n_1}} - j_2
x^{w_{n_2}})) = |\Delta_i| \cdot \|p\|^2.
$$
Consequently
\begin{eqnarray}
& & \left| \E(T_i^2|\F_{i-1}) - |\Delta_i| \cdot \|p\|^2 \right| \nonumber\\
& \leq & \sum_{n_1,n_2 \in \Delta_i} ~\sum_{j_1,j_2=1}^d \left| \E \left(\cos (2
\pi \xi (j_1 x^{w_{n_1}} + j_2 x^{w_{n_2}})) \Big| \F_{i-1} \right) \right|
\label{var1}\\
& & + \underbrace{\sum_{n_1,n_2 \in \Delta_i} ~\sum_{j_1,j_2=1}^d}_{(j_1,n_1)
\neq (j_2,n_2)} \left| \E \left(\cos (2 \pi \xi (j_1 x^{w_{n_1}} - j_2
x^{w_{n_2}})) \Big| \F_{i-1} \right) \right|. \label{var2}
\end{eqnarray}
Let $[\alpha,\beta)$ be any atom of $\F_{i-1}$. Using ~\eqref{igl},~\eqref{lengex2} and
Lemma~\ref{lemma4}, we see that for any function from~\eqref{var1}
\begin{eqnarray*}
& & \frac{1}{(\beta-\alpha)} \int_\alpha^\beta \cos (2 \pi \xi (j_1 x^{w_{n_1}} +
j_2 x^{w_{n_2}}))~dx \\
& \leq & \frac{2 i^6 w_{\max(i-1)} \alpha^{w_{\max(i-1)}}}{\xi w_{\min(i)}
\alpha^{w_{\min(i)}}}  \\
& \ll & i^6 A^{-i}.
\end{eqnarray*}
Thus for any function from~\eqref{var1} we have
$$
\left| \E \left(\cos (2 \pi (j_1 x^{w_{n_1}} + j_2 x^{w_{n_2}})) \Big| \F_{i-1} \right) \right|
\ll i^6 A^{-i},
$$
and consequently the whole double sum in~\eqref{var1} is bounded by
\begin{eqnarray} 
& & |\Delta_i|^2 d^2 i^6 A^{-i} \nonumber\\
& \ll &  i^{14} A^{-i}. \label{at1}
\end{eqnarray}
Now consider any function of the form $\cos (2 \pi \xi (j_1 x^{w_{n_1}} - j_2 x^{w_{n_2}}))$
from~\eqref{var2}. Then, using Lemma~\ref{lemma5} with $\eta = i^{-12}$ we know
that there exist three intervals $I_1,I_2$ and $I_3$ of total measure at least
$1 - 2 B \eta=1-2Bi^{-12}$, such that for any atom $[\alpha,\beta)$ of $\F_{i-1}$
which is completely contained in one of the intervals $I_1,I_2$ or $I_3$ we have
\begin{eqnarray}
& & \frac{1}{(\beta-\alpha)} \int_\alpha^\beta \cos (2 \pi \xi (j_1 x^{w_{n_1}} -
j_2 x^{w_{n_2}}))~dx \nonumber\\
& \leq & \frac{2 i^6 w_{\max(i-1)} \alpha^{w_{\max(i-1)}}}{\xi \eta
w_{\min(i)} \alpha^{w_{\min(i)}}}  \nonumber\\
& \ll & i^{18} A^{-i}. \label{atunter}
\end{eqnarray}
Since by~\eqref{lengex} the length of any atom of $\F_{i-1}$ is at most 
$$
2^{-m(i)} \leq A^{-i},
$$
the total measure of those atoms of $\F_{i-1}$ which are \emph{not} completely
contained in one of the intervals $I_1,I_2,I_3$ is at most
$$
\p([A,B]) - \p(I_1 \cup I_2 \cup I_3) + 6 A^{-i} \ll i^{-12}.
$$
Combining this with~\eqref{atunter} we obtain for any function from~\eqref{var2}
$$
\left| \E \left(\cos (2 \pi (j_1 x^{w_{n_1}} - j_2 x^{w_{n_2}})) \Big| \F_{i-1}
\right) \right| \ll i^{18} A^{-i} + i^{-12} \ll i^{-12}.
$$
Thus the double sum in~\eqref{var2} is bounded by
\begin{equation} \label{at2}
|\Delta_i|^2 d^2 i^{-12} \ll i^{-4},
\end{equation}
and together with~\eqref{at1} we conclude that
\begin{equation} \label{ti}
\left| \E(T_i^2|\F_{i-1}) - |\Delta_i| \cdot \|p\|^2 \right| \ll i^{-4}.
\end{equation}
Consequently, setting  
$$
V_M = \sum_{i=1}^M \E (Y_i^2 | \F_{i-1}) \qquad \textrm{and} \qquad b_M =
\sum_{i=1}^M |\Delta_i| \cdot \|p\|^2,
$$
we obtain by~\eqref{tiyi} and~\eqref{ti} that
\begin{eqnarray}
| V_M - b_M | & \leq & \sum_{i=1}^M \left( \left| \E(Y_i^2|\F_{i-1}) -
\E(T_i^2|\F_{i-1}) \right| + \left| \E(T_i^2|\F_{i-1}) - |\Delta_i| \cdot
\|p\|^2 \right| \right) \nonumber\\
& \ll & \sum_{i=1}^M \left(i^2 + i^{-4}\right) \nonumber\\
& \ll & M^{3}. \label{vm}
\end{eqnarray}
To prove Lemma~\ref{lemma2}, we observe that by Lemma~\ref{heyde} we have
\begin{equation}
\sup_{t \in \R} \left| \p \left( \frac{Y_1 + \dots + Y_M}{\sqrt{b_M}} <
t \right) - \Phi(t) \right| \leq K \left( \frac{\sum_{i=1}^M
\mathbb{E}Y_i^4 + \mathbb{E} \left( (V_M - b_M)^2 \right)}{b_M^2}
\right)^{1/5}, \label{heydeq}
\end{equation}
where $K$ is an absolute constant and $\Phi(t)$ denotes the standard normal
distribution function.\\

Using Lemma~\ref{lemma1} with $\delta = 2 (\log |\Delta_i|) (\log \log
|\Delta_i|)^{-1/2}$ we obtain
$$
\p \left( |T_i| \geq 2 c_A \sqrt{|\Delta_i|} \log |\Delta_i| \right) \ll
e^{-\log (i^4)} + \frac{1}{i^{64}} \ll i^{-64}.
$$
Consequently, since $|T_i| \ll |\Delta_i|$, we have
\begin{equation} \label{ti4}
\E (T_i^4) \ll \left(\sqrt{|\Delta_i|} \log |\Delta_i|\right)^4 + i^{-64}
|\Delta_i|^4 \ll i^8 (\log i)^4,
\end{equation}
and, by~\eqref{tiyi}, we obtain
\begin{equation} \label{tiyi2}
\E(Y_i^4) \ll i^8 (\log i)^4
\end{equation} 
and
\begin{equation}
\sum_{i=1}^M \E(Y_i^4) \ll \sum_{i=1}^M i^8 (\log i)^4 \ll M^9 (\log M)^4.
\end{equation}
On the other hand, 
\begin{equation} \label{bm}
b_M \gg \sum_{i=1}^M |\Delta_i| \gg M^5.
\end{equation}
Combining~\eqref{vm},~\eqref{heydeq},~\eqref{ti4} and~\eqref{bm} we get
$$
\sup_t \left| \p \left( \frac{\sum_{i=1}^M Y_i}{\sqrt{b_M}} <
t \right) - \Phi(t) \right| \ll \left( \frac{M^9 (\log M)^4 + M^6}{M^{10}}
\right)^{1/5} \ll \frac{\log M}{M^{1/5}}.
$$
By~\eqref{yzest} we have
\begin{equation} \label{fra1}
\left| \sum_{i=1}^M \left(T_i - Y_i\right) \right| \ll \sum_{i=1}^M i^{-2} \ll 1.
\end{equation}
Furthermore, we have
\begin{equation} \label{fra2}
\sqrt{N} \|p\| - \sqrt{b_M} \leq \frac{N \|p\|^2 - b_M}{\sqrt{N}\|p\| + \sqrt{b_M}} \ll \frac{\sum_{i=1}^M i \|p\|^2}{M^{5/2}} \ll M^{-1/2}, 
\end{equation}
and
\begin{equation} \label{fra3}
\left| T_1' + \dots T_M' \right| \ll \sum_{i=1}^M |\Delta_i'| \ll M^2.
\end{equation}
Note that
$$
\sum_{i=1}^M Y_i = \sum_{n=1}^N p(\xi x^{w_n}) + \sum_{i=1}^M \left(T_i - Y_i\right) + \sum_{i=1}^M T_i',
$$
and consequently
$$
\sum_{n=1}^N p(\xi x^{w_n}) < \|p\| \sqrt{N} t
$$
is equivalent to
$$
\frac{\sum_{i=1}^M Y_i}{\sqrt{b_M}} <  \underbrace{t \left( \frac{\|p\|\sqrt{N}}{\sqrt{b_M}} - \frac{\sum_{i=1}^M \left(T_i - Y_i\right)}{t\sqrt{b_M}} - \frac{\sum_{i=1}^M T_i'}{t\sqrt{b_M}} \right)}_{=: \hat{t}}.
$$
Using~\eqref{fra1},~\eqref{fra2} and~\eqref{fra3} we get
$$
\hat{t} = t \left(1 + \mathcal{O} \left(M^{-3} + M^{-1/2} t^{-1} \right)\right),
$$
and consequently
$$
\left| \Phi(t) - \Phi\left(\hat{t}\right) \right| \ll M^{-1/2}.
$$
Thus we finally get
$$
\sup_{t \in \R} \left| \p \left( \frac{\sum_{n=1}^N p(\xi x^{w_n})}{\|p\| \sqrt{N}} <
t \right) - \Phi(t) \right| \ll \frac{\log M}{M^{1/5}} \ll \frac{\log N}{N^{1/25}},
$$
which proves Lemma~\ref{lemma2}.\\

To prove Lemma~\ref{lemma2lil}, note that by~\eqref{vm},~\eqref{tiyi2} and~\eqref{bm} we have
$$
\frac{V_M}{b_M} \to 1
$$
and
$$
\sum_{M=1}^\infty \frac{(\log b_M)^{10}}{b_M^2} \E Y_M^4 \ll \sum_{M=1}^\infty \frac{(\log M)^{10}}{M^{10}} M^8 (\log M)^4 < \infty. 
$$
Thus by Lemma~\ref{lstr} we have
$$
\limsup_{M \to \infty} \frac{\sum_{i=1}^M Y_i}{\sqrt{2 b_M \log \log b_M}} = 1 \qquad \textup{a.s.}
$$
Since $\sum_{i=1}^M T_i \ll M^2 \ll \sqrt{b_M}$ and 
$$
\frac{b_M}{\|p\| \sum_{i=1}^M (|\Delta_i| + |\Delta_i'|)} \to 1 \qquad \textrm{as $M \to \infty$},
$$
we obtain, using~\eqref{yzest}, that
$$
\limsup_{M \to \infty} \frac{\sum_{i=1}^M (T_i + T_i')}{\sqrt{\sum_{i=1}^M (|\Delta_i| + |\Delta_i'|) \log \log \left( \sum_{i=1}^M (|\Delta_i| + |\Delta_i'|) \right)}} = \sqrt{2} \|p\| \qquad \textup{a.e.}
$$
which is Lemma~\ref{lemma2lil}.

\section{The law of the iterated logarithm for trigonometric polynomials}
\label{sectp}

In the present section we will prove the exact law of the iterated logarithm for
trigonometric polynomials. 

\begin{lemma} \label{lemmaup}
Let $p$ be a trigonometric polynomial. Then for almost all $x \in [A,B]$
$$
\limsup_{N \to \infty} \frac{\left|\sum_{n=1}^N p(\xi x^{s_n})\right|}{\sqrt{N
\log \log N}} = \sqrt{2} \|p\|.
$$
\end{lemma}

\emph{Proof of Lemma~\ref{lemmaup}:~} Choose $\theta>1$ (``small''). For any $k \geq 1$, let $N_k$ denote the smallest number of the form 
\begin{equation} \label{form1}
\sum_{i=1}^M (i^4+i) \qquad \textrm{for some $M$}
\end{equation}
which satisfies 
\begin{equation} \label{nktk}
N_k \geq \theta^k.
\end{equation}
Then by Lemma~\ref{lemma2lil} we have
\begin{equation} \label{lower}
\limsup_{k \to \infty} \frac{\left| \sum_{n=1}^{N_k} p(\xi x^{s_n}) \right|}{\sqrt{N_k \log \log N_k}} = \sqrt{2} \|p\| \qquad \textup{a.e.}
\end{equation}
Since the sequence of numbers of the form~\eqref{form1} grows polynomially, we have
\begin{equation} \label{nksize}
\frac{N_{k+1}}{N_k} \to \theta \quad \textrm{and} \quad \frac{N_{k}}{\theta^k} \to 1 \qquad \textrm{as $k \to \infty$}.
\end{equation}
Let
\begin{eqnarray*}
V_k & = & \left( x \in [A,B]:~\max_{N_k < M \leq N_{k+1}} \left| \sum_{n=N_k+1}^M p(\xi x^{s_n}) \right| \right. \\
& & \qquad\qquad   \geq 238 c_A \sqrt{(N_{k+1}-N_k) \log \log (N_{k+1}-N_k)} \Bigg), \quad k \geq 1,
\end{eqnarray*}
where $c_A$ is the constant from Lemma~\ref{lemma1}. Using Corollary~\ref{co3} with $\gamma=2 \|r\|^{-1/4}$ we obtain
$$
\p (V_k) \ll \exp(-2 \log \log (N_{k+1}-N_k)) \ll k^{-2}.
$$
Thus
$$
\sum_{k=1}^\infty \p (V_k) < \infty,
$$
which implies that with probability one only finitely many events $V_k$ occur. Hence, by~\eqref{lower} and~\eqref{nksize}, for almost all $x \in [A,B]$ we have
$$
\sqrt{2}\| p\| \leq \limsup_{N \to \infty} \frac{\left| \sum_{n=1}^N p (\xi x^{s_n}) \right|}{\sqrt{N \log \log N}} \leq \sqrt{2} \|p\| + 238 c_A \sqrt{1-\theta^{-1}}.
$$
Since $\theta>1$ can be chosen arbitrarily close to 1, this proves Lemma~\ref{lemmaup}.

\section{The law of the iterated logarithm for functions of bounded variation}
\label{sectf}

\begin{lemma} \label{lemmabv}
For any function $f$ satisfying~\eqref{f} we have for almost all $x \in [A,B]$ we have
$$
\limsup_{N \to \infty} \frac{\left| \sum_{k=1}^N f(\xi x^{s_n}) \right|}{\sqrt{N \log \log N}}
= \sqrt{2} \|f\|
$$
\end{lemma}

Choose $d \geq 1$, and split the Fourier series of $f$ into a trigonometric
polynomial of degree $d$ (which will be denoted by $p$) and a remainder function
$r$. Then, since
$$
p(x) - r(x) \leq f(x) \leq p(x) + r(x),
$$
by Lemma~\ref{lemmar} and Lemma~\ref{lemmaup} for almost all $x \in [A,B]$ we have
$$
\sqrt{2} \|p\| - 238 c_A d^{-1/8} \leq \limsup_{N \to \infty} \frac{\left| \sum_{k=1}^N f(\xi
x^{s_n}) \right| }{\sqrt{N \log \log N}} \leq \sqrt{2} \|p\| + 238 c_A d^{-1/8}.
$$
By Parseval's identity we have $\|p\| \to \|f\|$ as $d \to \infty$. Since $d$ can be chosen
arbitrarily, we obtain
$$
\limsup_{N \to \infty} \frac{\left| \sum_{k=1}^N f(\xi x^{s_n}) \right|}{\sqrt{N \log \log N}}
= \sqrt{2} \|f\|
$$for almost all $x \in [A,B]$, which proves the lemma.

\section{Proof of Theorem~\ref{th1}} \label{sectth1}

First we note that
$$
\sup_{0 \leq a < b \leq 1} \left\|\mathbf{I}_{[a,b)} \right\| = \|
\mathbf{I}_{[0,1/2)} \| = \frac{1}{2}.
$$
Thus by Lemma~\ref{lemmabv} for any $R \geq 1$
\begin{equation*}
\limsup_{N \to \infty} \frac{\sqrt{N} D_N^{(\geq 2^{-R})}(\{\xi
x^{s_n}\})}{\sqrt{\log \log N}} \geq \limsup_{N \to \infty} \frac{\left|
\sum_{k=1}^N \mathbf{I}_{[0,1/2)} (\xi x^{s_n}) \right|}{\sqrt{\log \log N}} =
\frac{1}{\sqrt{2}}
\end{equation*}
for almost all $x \in [A,B]$. Together with Lemma~\ref{lemmasmint} and
\eqref{discrepancies} this proves
$$
\frac{1}{\sqrt{2}} \leq \limsup_{N \to \infty} \frac{\sqrt{N} D_N^*(\{\xi
x^n\})}{\sqrt{\log \log N}} \leq \limsup_{N \to \infty} \frac{\sqrt{N} D_N(\{\xi
x^n\})}{\sqrt{\log \log N}} \leq \frac{1}{\sqrt{2}} + 3 \cdot 10^7 c_A R^{-1}
$$
for almost all $x \in [A,B]$. Since $R$ can be
chosen arbitrarily large, this proves Theorem~\ref{th1}.

\section{Proof of Theorem~\ref{th2}} \label{sectth2}

In this section we give a sketch of the proof of Theorem~\ref{th2}. This theorem is in large parts a consequence of Lemma~\ref{lemma2} from Section~\ref{sectmar}. By the assumptions made at the beginning of Section~\ref{pre} we
consider an interval $[A,B] \subset (1,\infty)$ which is of length 1; however, the proof remains true for intervals of arbitrary (positive) length in exactly the same way.\\

Let $\ve > 0$ be given. Let $N \geq 1$ also be given, and set 
$$
\hat{N} = \max\{ n \leq N:~\textrm{$n$ is of the form~\eqref{Nform}}\}.
$$
Then
\begin{equation} \label{nhut}
N - \hat{N} \ll N^{4/5}.
\end{equation}
Assume again for simplicity that $f$ is an even function, set $J = \lceil \ve^{-1} \rceil$ and
$$
p(x) = \sum_{j=1}^J a_j \cos 2 \pi j x, \qquad r(x) = \sum_{j=J+1}^\infty a_j \cos 2 \pi j x.
$$
Using the methods from Section~\ref{expi} we can prove that
\begin{equation} \label{rt}
\int_A^B \left( \sum_{n=1}^N r(\xi x^{s_n}) \right)^2~dx \ll N J^{-1}.
\end{equation}
Similarly, using~\eqref{nhut}, we can prove that
\begin{equation} \label{rt2}
\int_A^B \left( \sum_{n=\hat{N}+1}^N f(\xi x^{s_n}) \right)^2~dx \ll N^{4/5}.
\end{equation}
For the distribution of the normalized sum
\begin{equation} \label{rt4}
\frac{\sum_{n=1}^{\hat{N}} p(\xi x^{s_n})}{\|p\| \sqrt{\hat{N}}}
\end{equation}
we have the approximation result from Section~\ref{sectmar}. We further have
\begin{equation} \label{rt3}
\sum_{n=1}^{N} f(\xi x^{s_n}) = \sum_{n=1}^{\hat{N}} p(\xi x^{s_n}) + \sum_{n=\hat{N}+1}^N  p(\xi x^{s_n}) + \sum_{n=1}^N  r(\xi x^{s_n}).
\end{equation}
Now inequalities~\eqref{rt} and~\eqref{rt2} tell us that the last two sums on the right-hand side of~\eqref{rt3} are ``small'' with large probability in comparison to the normalizing factor $\sqrt{N}$. Note further that $\|p\| \to \|f\|$ as $\ve \to 0$. Thus the fact that the distribution of~\eqref{rt4} is close to the normal distribution tells us that also the distribution of 
$$
\frac{\sum_{n=1}^{N} f(\xi x^{s_n})}{\|f\| \sqrt{N}}
$$
is close to the normal distribution, provided $\ve$ is sufficiently small. Arguing as at the end of Section~\ref{sectmar}, all these results are sufficient to obtain
$$
\left| \p \left( x \in [A,B]:~\frac{\sum_{n=1}^N f(\xi x^{s_n})}{\|f\| \sqrt{N}} < t \right) - \Phi(t) \right| \ll \ve 
$$
for sufficiently large $N$, for all $t \in \R$.

\section*{Addendum}

I thank Katusi Fukuyama for the following argument, which illuminates the relation between Theorem~\ref{th1} and the corresponding results for the discrepancy of lacunary series. Let $x>1$. Then, according to~\eqref{fukuex}, we have
\begin{equation} \label{fu}
\frac{D_N(\{\xi x^n\})}{\sqrt{N \log \log N}} = \frac{1}{\sqrt{2}}
\qquad \textup{a.e. $\xi$,}
\end{equation}
if $x$ is a number for which $x^r \not\in \Q$ for all integers $r \geq 1$. Since the exceptional set of such $x$ is countable (and consequently has vanishing Lebesgue measure), this implies that 
$$
\frac{D_N(\{\xi x^n\})}{\sqrt{N \log \log N}} = \frac{1}{\sqrt{2}}
\qquad \textup{a.e. $\xi$, \quad for a.e. $x>1$.}
$$
Applying Fubini's theorem we obtain
\begin{equation} \label{wea}
\frac{D_N(\{\xi x^n\})}{\sqrt{N \log \log N}} = \frac{1}{\sqrt{2}}
\qquad \textup{a.e. $x>1$, \quad for a.e. $\xi$.}
\end{equation}
On the other hand, Corollary~\ref{co1} can be written as
$$
\frac{D_N(\{\xi x^n\})}{\sqrt{N \log \log N}} = \frac{1}{\sqrt{2}}
\qquad \textup{a.e. $x>1$, \quad for \emph{all} $\xi>0$.}
$$
Thus~\eqref{wea} is a significantly weaker version of Corollary~\ref{co1}, which, however, gives a plausible explanation of the appearance of the constant $1/\sqrt{2}$ in Corollary~\ref{co1}. A result similar to~\eqref{fu}, with the sequence $x^n$ replaced by $x^{s_n}$ for increasing $(s_n)_{n \geq 1}$, has been proved in~\cite{fukuhi}; applying again Fubini's theorem one can conclude that 
$$
\frac{D_N(\{\xi x^{s_n}\})}{\sqrt{N \log \log N}} = \frac{1}{\sqrt{2}}
\qquad \textup{a.e. $x>1$, \quad for a.e. $\xi$.}
$$
which is a weaker version of Theorem~\ref{th1}.


\end{document}